\documentclass[11pt,reqno]{amsart}

\usepackage[a4paper,margin=1in]{geometry}
\usepackage{amsmath,amssymb,amsthm,mathtools}
\usepackage[expansion=false]{microtype}
\usepackage[colorlinks=true, linkcolor=blue,%
citecolor=red, urlcolor=blue]{hyperref}

\usepackage{mathrsfs}
\allowdisplaybreaks
\newcommand{\Z}{\mathbb Z}
\newcommand{\Zge}{\mathbb Z_{\ge 0}}

\newcommand{\wt}{\operatorname{wt}}
\newcommand{\hw}{\operatorname{hw}}

\newcommand{\HH}{\mathcal H}

\numberwithin{table}{section}
\numberwithin{equation}{section}

\theoremstyle{plain}
\newtheorem{theorem}{Theorem}[section]
\newtheorem{lemma}[theorem]{Lemma}
\newtheorem{proposition}[theorem]{Proposition}
\newtheorem{corollary}[theorem]{Corollary}

\theoremstyle{definition}
\newtheorem{definition}[theorem]{Definition}
\newtheorem{remark}[theorem]{Remark}
\newtheorem{example}[theorem]{Example}
\title[Closed energy formulas in types $D_4^{(3)}$ and $G_2^{(1)}$]
{Closed energy formulas on tensor squares of perfect crystals of types $D_4^{(3)}$ and $G_2^{(1)}$}

\author{Shaolong Han}
\address{Beijing International Center for Mathematical Research, Peking University, Beijing 100871, P.R. China}
\email{hanshaolong@bicmr.pku.edu.cn}
\date{}

\begin{document}
\begin{abstract}
For each integer $l\ge1$, we give closed formulas for the local
energies on tensor squares of the level-$l$ perfect crystals of types $D_4^{(3)}$ and $G_2^{(1)}$, normalized to vanish at the all-zero coordinate tensors.
The $D_4^{(3)}$ energy is the maximum of three defects, equivalently of nine linear forms in parameters of classical highest-weight elements.
Eliminating these parameters yields the maximum of nine piecewise-linear functions of the twelve input coordinates.
Tensor-compatible similarity identifies the $G_2^{(1)}$ energy
at level $l$ with one third of the pullback of the $D_4^{(3)}$
energy at level $3l$, yielding the same coordinate formula.
The formula has no explicit level dependence, respects natural
coordinate-set inclusions between levels, and allows direct evaluation without the combinatorial $R$-matrix or crystal-graph traversal.
\end{abstract}

\maketitle

\section{Introduction}

Perfect crystals provide semi-infinite path realizations of highest-weight crystals for quantum affine algebras \cite{KMN2}.
Weighted sums of local energies determine the affine grading in the corresponding character formulas \cite{KMN1}.
On the tensor square of a perfect crystal, the local energy is constant along classical edges and changes along color-\(0\) edges according to the tensor factor selected by the Kashiwara operator.
Normalization determines this function uniquely, but its defining recursion does not make its coordinate dependence explicit.
The aim of this paper is to give closed coordinate formulas for the local energies on homogeneous tensor squares of the perfect-crystal families of types \(D_4^{(3)}\) and \(G_2^{(1)}\), uniformly at every positive level.

The underlying crystals are known at arbitrary level.
Kashiwara--Misra--Okado--Yamada constructed the type-\(D_4^{(3)}\) family \((B_l)_{l\ge1}\) in coordinates \cite{KMOY}.
Yamane constructed the type-\(G_2^{(1)}\) family \cite{Yamane}, and Misra--Mohamad--Okado subsequently made its color-\(0\) action explicit through similarity \cite{MMO}.

The energy functions have also been studied in several settings.
At level one, Benkart--Frenkel--Kang--Lee determined the energy for their uniform perfect-crystal construction in all affine types \cite{BFKL}.
For the two types considered here, energy tables describe the adjacent-column conditions in reduced Young walls \cite{FHKS}.
Dombos used the normalized type-\(D_4^{(3)}\) matrix to derive a grounded-partition identity \cite{Dombos}.

At arbitrary level, Yamada gave an insertion algorithm for the combinatorial \(R\)-matrix and energy on \(B_l\otimes B_1\) in type \(D_4^{(3)}\) \cite{YamadaSCA}, and Misra--Okado--Wilson treated the corresponding type-\(G_2^{(1)}\) products \cite{MOW}.
These products describe carrier--site interactions in soliton cellular automata.
Scrimshaw \cite{Scrimshaw2016} identified cocharge with total energy on arbitrary tensor products of the type-\(D_4^{(3)}\) crystals \(B_s\), denoted \(B^{1,s}\) in that work.
Our question concerns the direct evaluation of the homogeneous local energy in the twelve coordinates of its two input factors.

Closed homogeneous tensor-square energy formulas are known for seven classical affine families \cite{HanClassical}.
The two exceptional types considered here are outside that scope.
Moreover, formulas for \(B_l\otimes B_1\) do not directly give formulas for \(B_l\otimes B_l\): when \(l>1\), additional coordinate configurations in the second factor affect the string capacities and the thresholds for switching tensor factors.
To our knowledge, neither family's normalized homogeneous energy has previously been expressed uniformly at every level in the twelve input coordinates.

We normalize both energies to vanish at the corresponding all-zero coordinate tensors.
For type \(D_4^{(3)}\), six maximal classical raising strings determine a highest element with six parameters.
The energy is the maximum of three defect functions, equivalently the maximum of nine linear forms in these parameters (Theorem~\ref{thm:energy} and Proposition~\ref{prop:nine-branch}).
Evaluating the raising strings explicitly eliminates the highest parameters and gives the maximum of nine piecewise-linear functions of the twelve input coordinates (Theorem~\ref{thm:twelve-coordinate-energy}).

The proof turns the connectedness path of \cite[Section~7.1]{KMOY} into an energy evaluation.
Its first affine segment starts at a classical lowest tensor, whereas our parameters describe a highest tensor.
We connect these descriptions by computing a lowest normal form through six maximal lowering strings.
Defect coordinates then determine the number of affine moves on the left factor and its weighted coordinate sum at the endpoint, which yield the energy through the local-energy recursion.
A further difficulty arises when eliminating the highest parameters, since each raising stage depends on the changes made at earlier stages.
Explicit threshold functions resolve these dependencies and give coordinate formulas valid also at branch boundaries.

For type \(G_2^{(1)}\), we combine the similarity of \cite{MMO} with coordinate dilation into the level-\(3l\) type-\(D_4^{(3)}\) crystal.
To transfer the energy, we verify that every move in a similarity
block acts on the same tensor factor, including in the equality
case of the tensor product rule.
The normalized \(G_2^{(1)}\) energy is then one third of the pulled-back \(D_4^{(3)}\) energy (Theorem~\ref{thm:G-energy-transfer}).
This gives the highest-parameter formulas (Theorem~\ref{thm:G-closed-energy}), and positive homogeneity gives the same twelve-coordinate formula on the fractional coordinate lattice (Theorem~\ref{thm:G-direct-coordinate-energy}).
The transferred energy is integer-valued.

Thus a single piecewise-linear function describes the two energies on their respective coordinate domains.
The level enters only through these domains, and the formulas agree under coordinate-set inclusions between levels (Corollaries~\ref{cor:level-stability} and~\ref{cor:G-level-stability}).
They replace separate level-dependent energy tables by a fixed collection of coordinate operations, allowing direct evaluation without the combinatorial \(R\)-matrix or crystal-graph traversal.
The level-one examples recover the known energy matrices.

The paper is organized as follows.
Section~\ref{sec:crystal-preliminaries} fixes the crystal, tensor-product, and energy conventions.
Section~\ref{sec:perfect-crystal-models} records the affine data and coordinate models.
Section~\ref{sec:D43-closed-energy} derives the highest reduction, lowest normal form, defect dynamics, and energy formulas in type \(D_4^{(3)}\).
Section~\ref{sec:G21-similarity-transfer} proves tensor-product compatibility of similarity and transfers the formulas to type \(G_2^{(1)}\).

\section{Crystal-theoretic preliminaries}
\label{sec:crystal-preliminaries}

Let \(A=(a_{ij})_{i,j\in I}\) be a symmetrizable affine generalized
Cartan matrix with affine node \(0\), and let \(\mathfrak g\) be
the associated affine Kac--Moody algebra.
Write \(\alpha_i,h_i,\Lambda_i\) for its simple roots, simple
coroots, and affine fundamental weights, respectively, and
\(\delta_{\mathfrak g},c_{\mathfrak g}\) for its null root and
canonical central element.
Assume $\delta_{\mathfrak g}
=
\sum_{i\in I}a_i\alpha_i$ and $a_0=1$.

We use the lattices
\[
	P=
	\left(\bigoplus_{i\in I}\mathbb Z\Lambda_i\right)
	\oplus\mathbb Z\delta_{\mathfrak g},
	\quad
	P^\vee=\operatorname{Hom}_{\mathbb Z}(P,\mathbb Z),
	\quad
	P_{\mathrm{cl}}(\mathfrak g)
	:=
	P/\mathbb Z\delta_{\mathfrak g},
\]
with
\(\langle h_i,\alpha_j\rangle=a_{ij}\) and
\(\langle h_i,\Lambda_j\rangle=\delta_{ij}\).
We use the same symbols for roots and weights and their images in
\(P_{\mathrm{cl}}(\mathfrak g)\).
Since \(\langle c_{\mathfrak g},\delta_{\mathfrak g}\rangle=0\),
the level \(\langle c_{\mathfrak g},\lambda\rangle\) is well defined
for \(\lambda\in P_{\mathrm{cl}}(\mathfrak g)\).

For an indeterminate \(q\), let \(U_q(\mathfrak g)\) be the quantum
enveloping algebra over \(\mathbb Q(q)\), and let
\(U_q'(\mathfrak g)\) be its subalgebra without the degree operator.

\subsection{Crystals and tensor products}

Set \(P_{\mathrm{cl}}:=P_{\mathrm{cl}}(\mathfrak g)\) and
\(c:=c_{\mathfrak g}\).
On \(\mathbb Z\sqcup\{-\infty\}\), extend the usual order by
\(-\infty<n\) and set $(-\infty)\pm n=n+(-\infty)=-\infty$ $(n\in\mathbb Z)$.

Following \cite[Definition~2.1]{JayneMisra}, a
\(U_q'(\mathfrak g)\)-crystal is a set \(B\) with maps
\[
\wt:B\longrightarrow P_{\mathrm{cl}},
\qquad
\widetilde e_i,\widetilde f_i:
B\longrightarrow B\sqcup\{0\},
\qquad
\varepsilon_i,\varphi_i:
B\longrightarrow\mathbb Z\sqcup\{-\infty\}
\quad(i\in I),
\]
where \(0\notin B\) is a formal null element, satisfying the
following conditions for \(b,b'\in B\) and \(i\in I\):
\begin{enumerate}
	\item[\textup{(C1)}] $	\varphi_i(b)
	=
	\varepsilon_i(b)+\langle h_i,\wt(b)\rangle$.

	\item[\textup{(C2)}]
	If \(\widetilde e_i b\ne0\), then $(\wt,\varepsilon_i,\varphi_i)(\widetilde e_i b)
	=
	\bigl(
	\wt(b)+\alpha_i,\,
	\varepsilon_i(b)-1,\,
	\varphi_i(b)+1
	\bigr)$.

	\item[\textup{(C3)}]
	If \(\widetilde f_i b\ne0\), then $(\wt,\varepsilon_i,\varphi_i)(\widetilde f_i b)
	=
	\bigl(
	\wt(b)-\alpha_i,\,
	\varepsilon_i(b)+1,\,
	\varphi_i(b)-1
	\bigr)$.

	\item[\textup{(C4)}]
	$\widetilde f_i b=b'
	\iff
	\widetilde e_i b'=b$.

	\item[\textup{(C5)}]
	$\varphi_i(b)=-\infty
	\implies
	\widetilde e_i b=\widetilde f_i b=0$.
\end{enumerate}
We extend the operators by
\(\widetilde e_i0=\widetilde f_i0=0\).

All crystals considered below are finite and normal, meaning that
\begin{equation}\label{eq:normal-string-functions}
	\begin{aligned}
		\varepsilon_i(b)
		=
		\max\left\{
		m\in\mathbb Z_{\ge0}
		\mid
		\widetilde e_i^{\,m}b\ne0
		\right\},\quad
		\varphi_i(b)
		=
		\max \{
		m\in\mathbb Z_{\ge0}
		\mid
		\widetilde f_i^{\,m}b\ne0
		 \}
		\ (b\in B,\ i\in I),
	\end{aligned}
\end{equation}
see \cite[equation~(2.2.15)]{KMN1}.

For crystals \(B_1,B_2\), write
\(B_1\otimes B_2:=B_1\times B_2\), with elements
\(b_1\otimes b_2\).
Following \cite[equations~(2.2.16)--(2.2.19) and
Lemma~2.2.4(ii)]{KMN1}, its weight and string functions are
\begin{equation}\label{eq:tensor-string-functions}
	\begin{aligned}
		\wt(b_1\otimes b_2)
		&=
		\wt(b_1)+\wt(b_2),\\
		\varepsilon_i(b_1\otimes b_2)
		&=
		\max\left\{
		\varepsilon_i(b_1),
		\varepsilon_i(b_2)-\langle h_i,\wt(b_1)\rangle
		\right\},\\
		\varphi_i(b_1\otimes b_2)
		&=
		\max\left\{
		\varphi_i(b_2),
		\varphi_i(b_1)+\langle h_i,\wt(b_2)\rangle
		\right\}.
	\end{aligned}
\end{equation}
The Kashiwara operators are

\begin{equation}\label{eq:tensor-fe}
	\widetilde f_i(b_1\otimes b_2)
	=
	\begin{cases}
		\widetilde f_i b_1\otimes b_2,
		&\varphi_i(b_1)>\varepsilon_i(b_2),\\
		b_1\otimes\widetilde f_i b_2,
		&\varphi_i(b_1)\le\varepsilon_i(b_2),
	\end{cases}
	\quad
	\widetilde e_i(b_1\otimes b_2)
	=
	\begin{cases}
		\widetilde e_i b_1\otimes b_2,
		&\varphi_i(b_1)\ge\varepsilon_i(b_2),\\
		b_1\otimes\widetilde e_i b_2,
		&\varphi_i(b_1)<\varepsilon_i(b_2).
	\end{cases}
\end{equation}
A tensor containing \(0\) is identified with \(0\).

\subsection{Perfect crystals and local energy functions}

For a crystal \(B\), set $\varepsilon(b)
:=
\sum_{i\in I}\varepsilon_i(b)\Lambda_i$ and $\varphi(b)
:=
\sum_{i\in I}\varphi_i(b)\Lambda_i$ for $b\in B$,
and write
\(B_\lambda:=\{b\in B\mid\wt(b)=\lambda\}\).
Define
\[
\begin{aligned}
	P_{\mathrm{cl}}^+
	:=
	\left\{
	\lambda\in P_{\mathrm{cl}}
	\mid
	\langle h_i,\lambda\rangle\ge0
	\text{ for all }i\in I
	\right\},\quad
	(P_{\mathrm{cl}}^+)_l
	:=
	\left\{
	\lambda\in P_{\mathrm{cl}}^+
	\mid
	\langle c,\lambda\rangle=l
	\right\}
	\ (l\in\mathbb Z_{\ge0}).
\end{aligned}
\]
Fix \(l\in\mathbb Z_{>0}\), and set $B_{\min}
:=
\left\{
b\in B
\mid
\langle c,\varepsilon(b)\rangle=l
\right\}$.

\begin{definition}[{\cite[Definition~1.1.1]{KMN2}}]
	\label{def:perfect-crystal}
	A \(P_{\mathrm{cl}}\)-weighted crystal \(B\) is perfect of level
	\(l\) if:
	\begin{enumerate}
		\item[\textup{(P1)}]
		\(B\otimes B\) is connected.
		
		\item[\textup{(P2)}]
		There exists \(\lambda_0\in P_{\mathrm{cl}}\) such that $\wt(B)
		\subseteq
		\lambda_0+\sum_{i\ne0}\mathbb Z_{\le0}\alpha_i$ and $\#B_{\lambda_0}=1$.
	
		\item[\textup{(P3)}]
		There exists a finite-dimensional integrable
		\(P_{\mathrm{cl}}\)-weighted \(U_q'(\mathfrak g)\)-module
		with a crystal pseudobase \((L,B')\) such that
		\(B\cong B'/\{\pm1\}\).
		
		\item[\textup{(P4)}]
		\(\langle c,\varepsilon(b)\rangle\ge l\) for every \(b\in B\).
		
		\item[\textup{(P5)}]
		The maps
		\(\varepsilon,\varphi:
		B_{\min}\longrightarrow(P_{\mathrm{cl}}^+)_l\)
		are bijective.
	\end{enumerate}
\end{definition}

\begin{definition}\label{def:energy}
	Let \(B\) be a perfect crystal.
	A map \(H:B\otimes B\longrightarrow\mathbb Z\) is a local energy
	function if, for \(i\in I\) and \(X=b_1\otimes b_2\) with
	\(\widetilde f_iX\ne0\),
	\begin{equation}\label{eq:energy-def}
		H(\widetilde f_iX)-H(X)
		=
		\begin{cases}
			0,
			&i\ne0,\\
			-1,
			&i=0,\quad
			\varphi_0(b_1)>\varepsilon_0(b_2),\\
			1,
			&i=0,\quad
			\varphi_0(b_1)\le\varepsilon_0(b_2).
		\end{cases}
	\end{equation}
\end{definition}

\begin{proposition}\label{prop:perfect-energy-existence}
	Every perfect crystal \(B\) admits a local energy function,
	unique up to an additive constant.
\end{proposition}

\begin{proof}
	Let \(\lambda_0\) be as in \textup{(P2)}.
	The images of \(\alpha_i\), \(i\ne0\), are
	\(\mathbb Z\)-linearly independent in \(P_{\mathrm{cl}}\);
	hence \textup{(P2)} gives
	\begin{equation}\label{eq:unique-top-tensor-weight}
		(B\otimes B)_{2\lambda_0}
		=
		B_{\lambda_0}\otimes B_{\lambda_0},
		\qquad
		\#(B\otimes B)_{2\lambda_0}=1.
	\end{equation}
	
	Conditions \textup{(P1)}, \textup{(P3)}, and
	\eqref{eq:unique-top-tensor-weight} verify
	\cite[conditions~(3.4.2)--(3.4.3) and~(4.2.3)--(4.2.4)]{KMN1}.
	Thus \cite[Propositions~4.3.1 and~4.2.2]{KMN1} give an
	integer-valued energy function \(H\) on \(B\otimes B\)
	in the \(\widetilde e_i\)-convention.
	Reversing nonzero arrows using \eqref{eq:tensor-fe}
	yields \eqref{eq:energy-def}.
	
	For two local energy functions \(H,H'\),
	\eqref{eq:energy-def} makes \(H-H'\) constant along every edge.
	Connectedness in \textup{(P1)} therefore makes \(H-H'\)
	constant on \(B\otimes B\).
\end{proof}

\section{\texorpdfstring{The level-$l$ perfect crystals of types
		$D_4^{(3)}$ and $G_2^{(1)}$}
	{The level-l perfect crystals of types D4(3) and G2(1)}}
\label{sec:perfect-crystal-models}

\subsection{\texorpdfstring{Affine data for $D_4^{(3)}$ and $G_2^{(1)}$}
	{Affine data for D4(3) and G2(1)}}

Both types have index set \(I=\{0,1,2\}\).
Superscripts \(D,G\) distinguish their data and are omitted when
the type is fixed.
For \(t\in\mathbb R\), set \((t)_+:=\max\{t,0\}\).

For \(X\in\{D,G\}\), write
\(A^X=(a_{ij}^X)_{i,j\in I}\), where
\(a_{ij}^X=\langle h_i^X,\alpha_j^X\rangle\).
In the conventions of \cite[Section~2.1]{KMOY} and
\cite[Section~3.1.6]{Yamane},
\begin{equation}\label{eq:Cartan-matrices}
	A^D=
	\begin{pmatrix}
		2&-1&0\\
		-1&2&-3\\
		0&-1&2
	\end{pmatrix},
	\qquad
	A^G=
	\begin{pmatrix}
		2&-1&0\\
		-1&2&-1\\
		0&-3&2
	\end{pmatrix}.
\end{equation}
The null roots and canonical central elements are
\[
\begin{aligned}
	\delta^D
	&=\alpha_0^D+2\alpha_1^D+\alpha_2^D,
	&
	c^D
	&=h_0^D+2h_1^D+3h_2^D,\\
	\delta^G
	&=\alpha_0^G+2\alpha_1^G+3\alpha_2^G,
	&
	c^G
	&=h_0^G+2h_1^G+h_2^G.
\end{aligned}
\]
We choose the affine fundamental weights so that $\alpha_j^X
=
\sum_{i\in I}a_{ij}^X\Lambda_i^X
+\delta_{j0}\delta^X$ $(X\in\{D,G\},\ j\in I)$.

The finite nodes \(1,2\) form type \(G_2\), with
\(\alpha_1^D,\alpha_2^G\) short and
\(\alpha_2^D,\alpha_1^G\) long.
For either type, set $\overline{\Lambda}_1^X
:=
\Lambda_1^X-2\Lambda_0^X$,
a level-zero representative of the first finite fundamental weight.
Write \(B^{G_2}(\lambda)\) for the highest-weight
\(U_q(G_2)\)-crystal of dominant integral weight \(\lambda\),
with the inherited node labels.

Following \cite[Section~1]{MMO}, the similarity factors satisfy
\begin{equation}\label{eq:G-similarity-factors}
	(m_0,m_1,m_2)=(3,3,1),
	\quad
	a_{ij}^G=\frac{m_j}{m_i}a_{ij}^D
	\ (i,j\in I).
\end{equation}

\subsection{\texorpdfstring{The level-$l$ perfect crystal of type
		$D_4^{(3)}$}
	{The level-l perfect crystal of type D4(3)}}

We use the coordinate realization of
\cite[Sections~4.1--4.2]{KMOY}.
For
\(b=(x_1,x_2,x_3,\bar x_3,\bar x_2,\bar x_1)\in\mathbb R^6\),
set
\begin{equation}\label{eq:s-def}
	s(b)
	:=
	x_1+x_2+\frac{x_3+\bar x_3}{2}+\bar x_2+\bar x_1.
\end{equation}
For \(l\in\Z_{>0}\), define
\begin{equation}\label{eq:Bl-set}
	B_l
	:=
	\left\{
	b\in\Zge^6
	\ \middle|\
	x_3\equiv\bar x_3\! \! \!\!\pmod 2,\  s(b)\le l
	\right\}.
\end{equation}
Its classical restriction is
\begin{equation}\label{eq:G2-decomp}
	B_l
	\cong
	\bigoplus_{j=0}^l
	B^{G_2}(j\overline{\Lambda}_1).
\end{equation}

For \(b\in\mathbb R^6\), set
\begin{equation}\label{eq:z-def}
	z_1:=\bar x_1-x_1,
	\quad
	z_2:=\bar x_2-\bar x_3,
	\quad
	z_3:=x_3-x_2,
	\quad
	z_4:=\frac{\bar x_3-x_3}{2},
\end{equation}
and
\begin{equation}\label{eq:A-vector}
	\mathcal A(b)
	=
	\bigl(A_1(b),\ldots,A_6(b)\bigr)
	:=
	\bigl(
	0,\,
	z_1,\,
	z_1+z_2,\,
	z_1+z_2+3z_4,\,
	z_1+z_2+z_3+3z_4,\,
	2z_1+z_2+z_3+3z_4
	\bigr).
\end{equation}
Write $M(b):=\max_{1\le k\le6}A_k(b)$ and $D(b):=A_6(b)=2z_1+z_2+z_3+3z_4$.
By \eqref{eq:Bl-set}, we have $s(b),z_1,z_2,z_3,z_4,M(b),D(b)\in\mathbb Z$ $(b\in B_l)$.

For \(b\in B_l\) and \(1\le i\le6\), define
\begin{align}
	(F_i)
	&\iff
	A_i(b)=M(b)
	\ \text{ and }\
	A_j(b)<A_i(b)
	\quad\text{for all }j<i,
	\label{eq:F-max-char}\\
	(E_i)
	&\iff
	A_i(b)=M(b)
	\ \text{ and }\
	A_j(b)<A_i(b)
	\quad\text{for all }j>i.
	\label{eq:E-max-char}
\end{align}
These agree with \cite[Remark~4.1]{KMOY};
each \(b\in B_l\) satisfies exactly one \((F_i)\) and one \((E_i)\).

In the formulas below, omitted coordinates are unchanged.
Outputs outside \(B_l\) are interpreted as the formal null element
\(0\), distinct from the zero six-tuple in \(B_l\).

The color-\(0\) operators are
\begin{equation}\label{eq:f0}
	\widetilde f_0b
	=
	\begin{cases}
		(x_1+1,\ldots),
		&(F_1),\\
		(\ldots,x_3+1,\bar x_3+1,\ldots,\bar x_1-1),
		&(F_2),\\
		(\ldots,x_3+2,\ldots,\bar x_2-1,\ldots),
		&(F_3),\\
		(\ldots,x_2+1,\ldots,\bar x_3-2,\ldots),
		&(F_4),\\
		(x_1+1,\ldots,x_3-1,\bar x_3-1,\ldots),
		&(F_5),\\
		(\ldots,\bar x_1-1),
		&(F_6),
	\end{cases}
\end{equation}
\[
\widetilde e_0b
=
\begin{cases}
	(x_1-1,\ldots),
	&(E_1),\\
	(\ldots,x_3-1,\bar x_3-1,\ldots,\bar x_1+1),
	&(E_2),\\
	(\ldots,x_3-2,\ldots,\bar x_2+1,\ldots),
	&(E_3),\\
	(\ldots,x_2-1,\ldots,\bar x_3+2,\ldots),
	&(E_4),\\
	(x_1-1,\ldots,x_3+1,\bar x_3+1,\ldots),
	&(E_5),\\
	(\ldots,\bar x_1+1),
	&(E_6).
\end{cases}
\]

For the classical colors,
\begin{equation}\label{eq:e1}
	\widetilde e_1b
	=
	\begin{cases}
		(\ldots,\bar x_2+1,\bar x_1-1),
		&z_2\ge(-z_3)_+,\\
		(\ldots,x_3+1,\bar x_3-1,\ldots),
		&z_2<0\le z_3,\\
		(x_1+1,x_2-1,\ldots),
		&(z_2)_+<-z_3,
	\end{cases}
\end{equation}
\begin{equation}\label{eq:f1}
	\widetilde f_1b
	=
	\begin{cases}
		(x_1-1,x_2+1,\ldots),
		&(z_2)_+\le-z_3,\\
		(\ldots,x_3-1,\bar x_3+1,\ldots),
		&z_2\le0<z_3,\\
		(\ldots,\bar x_2-1,\bar x_1+1),
		&z_2>(-z_3)_+,
	\end{cases}
\end{equation}
\begin{equation}\label{eq:e2}
	\widetilde e_2b
	=
	\begin{cases}
		(\ldots,\bar x_3+2,\bar x_2-1,\ldots),
		&z_4\ge0,\\
		(\ldots,x_2+1,x_3-2,\ldots),
		&z_4<0,
	\end{cases}
\end{equation}
\begin{equation}\label{eq:f2}
	\widetilde f_2b
	=
	\begin{cases}
		(\ldots,x_2-1,x_3+2,\ldots),
		&z_4\le0,\\
		(\ldots,\bar x_3-2,\bar x_2+1,\ldots),
		&z_4>0.
	\end{cases}
\end{equation}

The string functions are
\begin{align}
	\varepsilon_1(b)
	&=
	\bar x_1+\bigl(-z_2+(-z_3)_+\bigr)_+,
	\quad
	\varphi_1(b)
	=
	x_1+\bigl(z_3+(z_2)_+\bigr)_+,
	\label{eq:epsphi1}\\
	\varepsilon_2(b)
	&=
	\bar x_2+(-z_4)_+,
	\quad
	\varphi_2(b)
	=
	x_2+(z_4)_+,
	\label{eq:epsphi2}\\
	\varepsilon_0(b)
	&=
	l-s(b)+M(b)-D(b),
	\quad
	\varphi_0(b)
	=
	l-s(b)+M(b).
	\label{eq:epsphi0}
\end{align}
The color-\(0\) formulas are
\cite[Proposition~4.4]{KMOY}.
By \eqref{eq:Bl-set} and \(M(b)\ge\max\{0,D(b)\}\),
all string functions take values in \(\Zge\).

The string functions and \cite[Lemma~7.1]{KMOY} give
\[
\langle c,\varepsilon(b)\rangle
=
\langle c,\varphi(b)\rangle
\ge l,
\quad
\langle c,\varepsilon(b)\rangle=l
\iff
z_1=z_2=z_3=z_4=0.
\]
Hence $(B_l)_{\min}
=
\left\{
(\alpha,\beta,\beta,\beta,\beta,\alpha)
\ \middle|\
\alpha,\beta\in\Zge,\ 
2\alpha+3\beta\le l
\right\}$.
For
\(b=(\alpha,\beta,\beta,\beta,\beta,\alpha)\in(B_l)_{\min}\), we have $\varepsilon(b)
=
\varphi(b)
=
(l-2\alpha-3\beta)\Lambda_0
+\alpha\Lambda_1+\beta\Lambda_2$.

By \cite[Theorem~1.2]{KMOY}, \(B_l\) is perfect of level \(l\).

\subsection{\texorpdfstring{The level-$l$ perfect crystal of type
		$G_2^{(1)}$}
	{The level-l perfect crystal of type G2(1)}}

We use the coordinate realization \(V_l\) of
\cite[Section~3]{MMO}, denoted here by \(B_l^G\).
For \(l\in\Z_{>0}\), define
\begin{equation}\label{eq:G-Bl-set}
	B_l^G
	:=
	\left\{
	b=(x_1,x_2,x_3,\bar x_3,\bar x_2,\bar x_1)
	\in\left(\tfrac13\Zge\right)^6
	\ \middle|\
	\begin{gathered}
		3x_3\equiv3\bar x_3\! \! \!\!\pmod 2,\\
		x_1,\bar x_1,x_2-x_3,\bar x_3-\bar x_2\in\mathbb Z,
		\  s(b)\le l
	\end{gathered}
	\right\}.
\end{equation}
By \cite[Proposition~2]{MMO}, its classical restriction is $B_l^G
\cong
\bigoplus_{j=0}^l
B^{G_2}(j\overline{\Lambda}_1)$.

For \(b\in B_l^G\), use \(z_1,\ldots,z_4\) from
\eqref{eq:z-def} and set
\[
\begin{aligned}
	\mathcal A_G(b)
	=
	\bigl(A_1^G(b),\ldots,A_6^G(b)\bigr)
	:=
	\mathcal A(b),\quad
	M_G(b):=M(b),
	\quad
	D_G(b):=D(b).
\end{aligned}
\]
The conditions in \eqref{eq:G-Bl-set} give
\[
\begin{aligned}
	s(b)
	=
	x_1+\bar x_1-z_3+z_2
	+\frac{3x_3+3\bar x_3}{2}
	\in\Zge,\quad
	z_1,z_2,z_3,3z_4,M_G(b),D_G(b)
	\in\mathbb Z.
\end{aligned}
\]
Define \((F_i)\) and \((E_i)\) on \(B_l^G\) by
\eqref{eq:F-max-char}--\eqref{eq:E-max-char}, using
\(A_k^G,M_G\).

We use the operators of \cite[Section~3, Proposition~1]{MMO}.
Omitted coordinates are unchanged, and outputs outside \(B_l^G\)
are interpreted as \(0\).
The color-\(0\) operators are
\begin{equation}\label{eq:G-e0}
	\widetilde e_0^G b
	=
	\begin{cases}
		(x_1-1,\ldots),
		&(E_1),\\
		(\ldots,x_3-1,\bar x_3-1,\ldots,\bar x_1+1),
		&(E_2),\\
		(\ldots,x_2-\tfrac23,x_3-\tfrac23,
		\bar x_3+\tfrac43,\bar x_2+\tfrac13,\ldots),
		&(E_3),\ z_4=-\tfrac13,\\
		(\ldots,x_2-\tfrac13,x_3-\tfrac43,
		\bar x_3+\tfrac23,\bar x_2+\tfrac23,\ldots),
		&(E_3),\ z_4=-\tfrac23,\\
		(\ldots,x_3-2,\ldots,\bar x_2+1,\ldots),
		&(E_3),\ z_4\le-1,\\
		(\ldots,x_2-1,\ldots,\bar x_3+2,\ldots),
		&(E_4),\\
		(x_1-1,\ldots,x_3+1,\bar x_3+1,\ldots),
		&(E_5),\\
		(\ldots,\bar x_1+1),
		&(E_6),
	\end{cases}
\end{equation}
\begin{equation}\label{eq:G-f0}
	\widetilde f_0^G b
	=
	\begin{cases}
		(x_1+1,\ldots),
		&(F_1),\\
		(\ldots,x_3+1,\bar x_3+1,\ldots,\bar x_1-1),
		&(F_2),\\
		(\ldots,x_3+2,\ldots,\bar x_2-1,\ldots),
		&(F_3),\\
		(\ldots,x_2+\tfrac13,x_3+\tfrac43,
		\bar x_3-\tfrac23,\bar x_2-\tfrac23,\ldots),
		&(F_4),\ z_4=\tfrac13,\\
		(\ldots,x_2+\tfrac23,x_3+\tfrac23,
		\bar x_3-\tfrac43,\bar x_2-\tfrac13,\ldots),
		&(F_4),\ z_4=\tfrac23,\\
		(\ldots,x_2+1,\ldots,\bar x_3-2,\ldots),
		&(F_4),\ z_4\ge1,\\
		(x_1+1,\ldots,x_3-1,\bar x_3-1,\ldots),
		&(F_5),\\
		(\ldots,\bar x_1-1),
		&(F_6).
	\end{cases}
\end{equation}

The right-hand sides of \eqref{eq:e1}--\eqref{eq:f1} define
\(\widetilde e_1^G\) and \(\widetilde f_1^G\), respectively.
For color \(2\),
\begin{equation}\label{eq:G-e2}
	\widetilde e_2^G b
	=
	\begin{cases}
		(\ldots,\bar x_3+\tfrac23,\bar x_2-\tfrac13,\ldots),
		&z_4\ge0,\\
		(\ldots,x_2+\tfrac13,x_3-\tfrac23,\ldots),
		&z_4<0,
	\end{cases}
\end{equation}
\begin{equation}\label{eq:G-f2}
	\widetilde f_2^G b
	=
	\begin{cases}
		(\ldots,x_2-\tfrac13,x_3+\tfrac23,\ldots),
		&z_4\le0,\\
		(\ldots,\bar x_3-\tfrac23,\bar x_2+\tfrac13,\ldots),
		&z_4>0.
	\end{cases}
\end{equation}

The string functions for colors \(0,1\) are given by
\eqref{eq:epsphi0} and \eqref{eq:epsphi1}, respectively,
with \(M,D\) replaced by \(M_G,D_G\).
For color \(2\),
\begin{equation}\label{eq:G-epsphi2}
	\varepsilon_2^G(b)
	=
	3\bigl(\bar x_2+(-z_4)_+\bigr),
	\quad
	\varphi_2^G(b)
	=
	3\bigl(x_2+(z_4)_+\bigr).
\end{equation}
These agree with
\cite[Theorem~2, equation~(3.3)]{MMO}.
By \eqref{eq:G-Bl-set} and
\(M_G(b)\ge\max\{0,D_G(b)\}\), we have
\(\varepsilon_i^G(b),\varphi_i^G(b)\in\Zge\) for \(i\in I\).

The string functions and \cite[Lemma~2]{MMO} give
\[
\langle c,\varepsilon(b)\rangle
=
\langle c,\varphi(b)\rangle
\ge l,
\quad
\langle c,\varepsilon(b)\rangle=l
\iff
z_1=z_2=z_3=z_4=0.
\]
Hence $(B_l^G)_{\min}
=
\left\{
(\alpha,\beta,\beta,\beta,\beta,\alpha)
\ \middle|\
\alpha\in\Zge,\ 
\beta\in\tfrac13\Zge,\ 
2\alpha+3\beta\le l
\right\}$.

For
\(b=(\alpha,\beta,\beta,\beta,\beta,\alpha)\in(B_l^G)_{\min}\), we have $\varepsilon(b)
=
\varphi(b)
=
(l-2\alpha-3\beta)\Lambda_0
+\alpha\Lambda_1+3\beta\Lambda_2$.

By \cite[Theorem~2]{MMO}, \(B_l^G\) is isomorphic to Yamane's
level-\(l\) crystal, which is perfect by
\cite[Theorem~2.2]{Yamane}.

\section{\texorpdfstring{Closed local energy in type $D_4^{(3)}$}
	{Closed local energy in type D4(3)}}
\label{sec:D43-closed-energy}

Fix \(l\in\Z_{>0}\), and set
\(\phi:=(0,0,0,0,0,0)\in B_l\).
By Proposition~\ref{prop:perfect-energy-existence}, there is a unique
local energy function
\(H_l:B_l\otimes B_l\longrightarrow\mathbb Z\)
satisfying
\begin{equation}\label{eq:H-normalization}
	H_l(\phi\otimes\phi)=0.
\end{equation}

\subsection{\texorpdfstring{Explicit $G_2$-highest reduction}
	{Explicit G2-highest reduction}}

By \eqref{eq:energy-def}, \(H_l\) is constant on each classical
connected component of \(B_l\otimes B_l\).
Proposition~\ref{prop:hw-highest} gives an explicit reduction to
its unique \(G_2\)-highest element.

\subsubsection{Closed formulas for the classical raising strings}

Let \(b=(x_1,x_2,x_3,\bar x_3,\bar x_2,\bar x_1)\in B_l\).
By \eqref{eq:epsphi1}, we have $\varepsilon_1(b)
=
\bar x_1+
\bigl(-z_3-(z_2)_+\bigr)_+
+(-z_2)_+$.
For \(m\in\Zge\) with \(m\le\varepsilon_1(b)\), put
\[
r_1:=\min\bigl\{m,\bigl(-z_3-(z_2)_+\bigr)_+\bigr\},
\quad
r_2:=\min\{m-r_1,(-z_2)_+\},
\quad
r_3:=m-r_1-r_2,
\]
and define
\begin{equation}\label{eq:E1m}
	E_1(b;m)
	:=
	\bigl(
	x_1+r_1,\,
	x_2-r_1,\,
	x_3+r_2,\,
	\bar x_3-r_2,\,
	\bar x_2+r_3,\,
	\bar x_1-r_3
	\bigr).
\end{equation}

For \(m\in\Zge\) with \(m\le\varepsilon_2(b)\), put $q_1:=\min\{m,(-z_4)_+\}$, $q_2:=m-q_1$,
and define
\begin{equation}\label{eq:E2m}
	E_2(b;m)
	:=
	\bigl(
	x_1,\,
	x_2+q_1,\,
	x_3-2q_1,\,
	\bar x_3+2q_2,\,
	\bar x_2-q_2,\,
	\bar x_1
	\bigr).
\end{equation}

\begin{lemma}\label{lem:closed-string}
	Let \(b\in B_l\), \(i\in\{1,2\}\), and
	\(m\in\Zge\) with \(m\le\varepsilon_i(b)\).
	Then \(E_i(b;m)=\widetilde e_i^{\,m}b\).
\end{lemma}

\begin{proof}
	By \eqref{eq:normal-string-functions}, the first \(m\) raising
	steps are nonzero.
	
	For \(i=1\), the third branch of \eqref{eq:e1} fixes \(z_2\)
	and increases \(z_3\) by \(1\) until \(z_3\ge-(z_2)_+\).
	The second branch then applies exactly when \(z_2<0\);
	in this case \(z_3\ge0\), and
	\[
	(z_2,z_3)\longmapsto(z_2+1,z_3+1)
	\]
	until \(z_2=0\).
	The first branch then fixes \(z_3\) and sends
	\[
	(z_2,\bar x_1)\longmapsto(z_2+1,\bar x_1-1);
	\]
	its condition persists until \(\bar x_1=0\).
	Thus the successive block lengths are $\bigl(-z_3-(z_2)_+\bigr)_+$, $	(-z_2)_+$ and $\bar x_1$.

	For \(i=2\), both branches of \eqref{eq:e2} send
	\(z_4\mapsto z_4+1\).
	Since \(z_4\in\mathbb Z\), the second branch has length
	\((-z_4)_+\).
	The first then has length \(\bar x_2\), decreasing \(\bar x_2\)
	by \(1\) at each step while keeping \(z_4\ge0\).
	
	Truncating these blocks after \(m\) steps gives
	\(r_1,r_2,r_3\) and \(q_1,q_2\), respectively.
	The coordinate updates in \eqref{eq:e1} and \eqref{eq:e2}
	therefore give \eqref{eq:E1m}--\eqref{eq:E2m}.
\end{proof}

\subsubsection{Highest reduction on the tensor square}

For \(i\in\{1,2\}\) and \(Y\) in either \(B_l\) or
\(B_l\otimes B_l\), set $\widetilde e_i^{\max}Y
:=
\widetilde e_i^{\,\varepsilon_i(Y)}Y$.
For \(X=b\otimes c\in B_l\otimes B_l\), define $\mathcal P_i(X)
:=
E_i\bigl(b;\varepsilon_i(b)\bigr)
\otimes
E_i\bigl(c;(\varepsilon_i(c)-\varphi_i(b))_+\bigr)$.
Normality gives $0\le(\varepsilon_i(c)-\varphi_i(b))_+
\le\varepsilon_i(c)$,
so \(\mathcal P_i\) is well defined by
Lemma~\ref{lem:closed-string}.

\begin{lemma}\label{lem:Pi-emax}
	For \(X\in B_l\otimes B_l\) and \(i\in\{1,2\}\),
	\(\mathcal P_i(X)=\widetilde e_i^{\max}X\).
\end{lemma}

\begin{proof}
	Write \(X=b\otimes c\), and put
	\(d:=(\varepsilon_i(c)-\varphi_i(b))_+\).
	For \(0\le k<d\),
	\[
	\varepsilon_i(\widetilde e_i^{\,k}c)
	=
	\varepsilon_i(c)-k
	>
	\varphi_i(b).
	\]
	Thus the first \(d\) moves act on the right factor.
	For \(0\le k<\varepsilon_i(b)\),
	\[
	\varepsilon_i(\widetilde e_i^{\,d}c)
	=
	\min\{\varepsilon_i(c),\varphi_i(b)\}
	\le
	\varphi_i(b)+k
	=
	\varphi_i(\widetilde e_i^{\,k}b),
	\]
	so the next \(\varepsilon_i(b)\) moves act on the left factor.
	
	It follows from \eqref{eq:tensor-string-functions} and \textup{(C1)} that $\varepsilon_i(X)
	=
	\varepsilon_i(b)
	+
	(\varepsilon_i(c)-\varphi_i(b))_+
	=
	\varepsilon_i(b)+d$.
	Lemma~\ref{lem:closed-string} therefore gives $\widetilde e_i^{\max}X
	=
	\widetilde e_i^{\,\varepsilon_i(b)}b
	\otimes
	\widetilde e_i^{\,d}c
	=
	\mathcal P_i(X)$.
\end{proof}

Put \((i_1,\ldots,i_6):=(1,2,1,2,1,2)\).
For \(X\in B_l\otimes B_l\), define
\begin{equation}\label{eq:HW-recursion}
	X^{(0)}:=X,
	\quad
	X^{(k)}
	:=
	\mathcal P_{i_k}\bigl(X^{(k-1)}\bigr)
	\ (1\le k\le6),
\end{equation}
and set
\begin{equation}\label{eq:hw-def}
	\hw(X):=X^{(6)}.
\end{equation}

For a dominant integral weight \(\lambda\) of the finite
\(G_2\)-subalgebra, let \(B(\pi_\lambda)\) denote the path crystal
generated by \(\pi_\lambda(t):=t\lambda\), \(0\le t\le1\),
under Littelmann's root operators.

\begin{proposition}\label{prop:hw-highest}
	For every \(X\in B_l\otimes B_l\), the element \(\hw(X)\) is the
	unique element in the classical connected component of \(X\)
	annihilated by both \(\widetilde e_1\) and \(\widetilde e_2\).
\end{proposition}

\begin{proof}
	By \eqref{eq:G2-decomp} and
	\cite[Theorems~2.4.2--2.4.3]{KMN1}, the classical component \(C\)
	containing \(X\) is isomorphic to \(B^{G_2}(\lambda)\) for some
	dominant integral weight \(\lambda\).
	By \cite[Theorem~4.1]{KashiwaraSimilarity}, identify \(C\) with
	\(B(\pi_\lambda)\), sending its unique highest element \(u_C\)
	to \(\pi_\lambda\).
	
	Apply
	\cite[Lemma~11.1(c), Theorem~11.3 and its proof]{LittelmannRO}
	with \(\pi_0=\pi_\lambda\), \(\gamma=\alpha_1^D\), and
	\(\alpha=\alpha_2^D\).
	By \eqref{eq:Cartan-matrices} and
	\cite[Lemma~11.1(c)]{LittelmannRO}, we have \(y_4=0\).
	The proof of \cite[Theorem~11.3]{LittelmannRO} constructs
	successive maximal raising strings from \(X\) to \(u_C\),
	in colors \(1,2,1,2,\ldots\), with exponents
	\(n_1,m_1,n_2,m_2,\ldots\), allowing \(n_1=0\).
	A nonempty seventh string would satisfy $1\le n_4\le3m_3y_4=0$.

	Thus \(u_C\) is reached within six stages.
	Appending zero exponents after reaching \(u_C\),
	Lemma~\ref{lem:Pi-emax} and
	\eqref{eq:HW-recursion}--\eqref{eq:hw-def} give $\hw(X)
	=
	X^{(6)}
	=
	\widetilde e_{i_6}^{\max}\cdots
	\widetilde e_{i_1}^{\max}X
	=
	u_C$.
\end{proof}

\begin{lemma}\label{lem:highest-shape}
	An element \(Y\in B_l\otimes B_l\) is annihilated by
	\(\widetilde e_1\) and \(\widetilde e_2\) if and only if there are
	unique \(a,u,v,w,d,k\in\Zge\) such that
	\begin{equation}\label{eq:highest-shape}
		Y
		=
		(a,0,0,0,0,0)
		\otimes
		(u,v,w,w+2d,0,k),
	\end{equation}
	\(a\le l\), \(u+v+w+d+k\le l\), and
	\begin{equation}\label{eq:highest-admissibility}
		k+w+2d+(v-w)_+\le a.
	\end{equation}
\end{lemma}

\begin{proof}
	Write \(Y=b\otimes c\).
	By \eqref{eq:normal-string-functions}--\eqref{eq:tensor-string-functions}  and \textup{(C1)},
	\[
	\widetilde e_iY=0
	\iff
	\varepsilon_i(b)=0
	\ \text{and}\
	\varepsilon_i(c)\le\varphi_i(b)
	\ (i=1,2).
	\]
	
	Write \(b=(x_1,x_2,x_3,\bar x_3,\bar x_2,\bar x_1)\).
	By \eqref{eq:epsphi1}--\eqref{eq:epsphi2},
	\(\varepsilon_2(b)=0\) implies $\bar x_2=0$, $x_3\le\bar x_3$ and $\varepsilon_1(b)
	=
	\bar x_1+\bar x_3+(x_2-x_3)_+$.
	Hence
	\[
	\varepsilon_1(b)=\varepsilon_2(b)=0
	\iff
	b=(a,0,0,0,0,0)
	\quad\text{for some }a\in\Zge,
	\]
	and then
	\((\varphi_1(b),\varphi_2(b))=(a,0)\).
	
	By \eqref{eq:epsphi2} and the parity condition in
	\eqref{eq:Bl-set},
	\[
	\varepsilon_2(c)=0
	\iff
	c=(u,v,w,w+2d,0,k)
	\quad\text{for some }u,v,w,d,k\in\Zge.
	\]
	For these factors, direct substitution gives $s(b)=a$, $s(c)=u+v+w+d+k$ and $\varepsilon_1(c)=k+w+2d+(v-w)_+$.

	The level bounds follow from \eqref{eq:Bl-set}, and the
	remaining condition
	\(\varepsilon_1(c)\le\varphi_1(b)=a\)
	is precisely \eqref{eq:highest-admissibility}.
	Uniqueness follows directly from the coordinates in
	\eqref{eq:highest-shape}.
\end{proof}

\subsection{Highest parameters and classical lowest normal form}

\subsubsection{Highest parameters and the energy candidate}

By Proposition~\ref{prop:hw-highest} and
Lemma~\ref{lem:highest-shape}, for every
\(X\in B_l\otimes B_l\) there are unique
\(a,u,v,w,d,k\in\Zge\) such that
\begin{equation}\label{eq:high-param}
	\hw(X)
	=
	(a,0,0,0,0,0)\otimes c^\sharp,
	\quad
	c^\sharp:=(u,v,w,w+2d,0,k).
\end{equation}
These parameters depend only on the classical connected component
of \(X\).

Set
\begin{equation}\label{eq:r-def}
	r:=s(c^\sharp)=u+v+w+d+k.
\end{equation}
Let \(D^\sharp:=D(c^\sharp)\) and \(M^\sharp:=M(c^\sharp)\).
By \eqref{eq:z-def}--\eqref{eq:A-vector}, and
\(k-u-w-2d\le k-u\), we have
\begin{align}
	D^\sharp
	&=2k-2u-v+d,
	\notag\\
	M^\sharp
	&=
	\max\{
	0,\,
	k-u,\,
	k-u-w+d,\,
	k-u-v+d,\,
	D^\sharp
	\}.
	\label{eq:M-sharp}
\end{align}
Set
\begin{equation}\label{eq:chi-sharp}
	\chi^\sharp
	:=
	2k+w+d+\min\{v,w,d\},
\end{equation}
and define
\begin{equation}\label{eq:UVW}
	U(X):=a+M^\sharp-r-D^\sharp,
	\quad
	V(X):=r-a,
	\quad
	W(X):=r-\chi^\sharp.
\end{equation}
The energy candidate
\(\HH:B_l\otimes B_l\longrightarrow\mathbb Z\) is
\begin{equation}\label{eq:H-three}
	\HH(X):=\max\{U(X),V(X),W(X)\}.
\end{equation}

\begin{proposition}
	\label{prop:nine-branch}
	For every \(X\in B_l\otimes B_l\),
	\begin{equation}\label{eq:H-nine}
		\begin{aligned}
			\HH(X)=\max\bigl\{&
			a+u-w-2d-3k,\,
			a-w-2d-2k,\,
			a-2w-d-2k,\\
			&
			a-v-w-d-2k,\,
			a-r,\,
			r-a,\,
			u-k,\,
			u+v-w-k,\,
			u+v-d-k
			\bigr\}.
		\end{aligned}
	\end{equation}
\end{proposition}

\begin{proof}
	Equations~\eqref{eq:r-def}--\eqref{eq:M-sharp} give $a-r-D^\sharp=a+u-w-2d-3k$.
	Substituting the expansion of \(M^\sharp\) from
	\eqref{eq:M-sharp} into \eqref{eq:UVW} gives the first five
	branches in \eqref{eq:H-nine}.
	Moreover, \(V(X)=r-a\), and
	\eqref{eq:r-def}, \eqref{eq:chi-sharp}--\eqref{eq:UVW} give
	\[
	\begin{aligned}
		W(X)
		=u+v-k-\min\{v,w,d\}
		=\max\{u-k,u+v-w-k,u+v-d-k\}.
	\end{aligned}
	\]
	Equation~\eqref{eq:H-three} now yields \eqref{eq:H-nine}.
\end{proof}

\subsubsection{The classical lowest normal form}

For \(X\in B_l\otimes B_l\), we compute the classical lowest
tensor in its component from the highest parameters
\eqref{eq:high-param}, following the connectedness path in
\cite[Section~7.1]{KMOY}.

Put $p:=(v-w)_+$, $q:=(w-v)_+$ and $\mu:=\min\{v,w\}$.
Then
\begin{equation}\label{eq:vw-pqmu}
	v=\mu+p,
	\quad
	w=\mu+q,
	\quad
	pq=0.
\end{equation}
By \eqref{eq:highest-admissibility}, we have
\begin{equation}\label{eq:h-def}
	h:=a-k-w-2d-p\in\Zge.
\end{equation}

For \(m\in\Zge\), write
\(\bar{\mathbf1}^{\,m}:=(0,0,0,0,0,m)\), and define
\begin{equation}\label{eq:bdown-def}
	b^\downarrow
	:=
	(k,0,2\mu+q+d,q+d,d+p,h).
\end{equation}
Its coordinates are nonnegative integers, with
\(2\mu+q+d\equiv q+d\pmod2\), and
\begin{equation}\label{eq:s-bdown}
	s(b^\downarrow)
	=
	k+\mu+q+2d+p+h
	=
	a.
\end{equation}
Since \(a,r\le l\) by Lemma~\ref{lem:highest-shape},
both \(b^\downarrow\) and \(\bar{\mathbf1}^{\,r}\) belong to \(B_l\).
Set $X^\downarrow
:=
b^\downarrow\otimes\bar{\mathbf1}^{\,r}$.

\begin{lemma}
	\label{lem:classical-lowest-normal-form}
	For every \(X\in B_l\otimes B_l\), the element \(X^\downarrow\)
	is the unique element in the classical component of \(X\)
	annihilated by \(\widetilde f_1\) and \(\widetilde f_2\).
\end{lemma}

\begin{proof}
	Put $\tau:=(\mu-q-d)_+$, $\sigma:=(q+d-\mu)_+$ and $\eta:=\min\{\mu,q+d\}$.
	Then
	\begin{equation}\label{eq:tau-sigma-eta-identities}
		\tau+\eta=\mu,
		\quad
		\sigma+\eta=q+d,
		\quad
		\tau\sigma=0.
	\end{equation}
	
	Set \(Y^{(0)}:=\hw(X)=L_0\otimes R_0\), and define
	\begin{equation}\label{eq:classical-lowering-recursion}
		Y^{(j)}
		:=
		\widetilde f_{i_j}^{
			\,\varphi_{i_j}(Y^{(j-1)})}
		Y^{(j-1)}
		=:L_j\otimes R_j
		\ (1\le j\le6).
	\end{equation}
	
	Let
	\(\mathsf I_1,\mathsf{II}_1,\mathsf{III}_1\) and
	\(\mathsf I_2,\mathsf{II}_2\) denote the branches of
	\eqref{eq:f1} and \eqref{eq:f2}, respectively, in their displayed
	order.
	For a one-factor string starting at
	\(b=(x_1,x_2,x_3,\bar x_3,\bar x_2,\bar x_1)\in B_l\),
	the branch conditions give the full branch words
	\[
	\begin{aligned}
		&\mathsf{III}_1^{(z_2-(-z_3)_+)_+}
		\mathsf{II}_1^{(z_3)_+}
		\mathsf I_1^{x_1}
		&&\text{for color }1,\\
		&\mathsf{II}_2^{(z_4)_+}
		\mathsf I_2^{x_2}
		&&\text{for color }2.
	\end{aligned}
	\]
	All exponents are evaluated at \(b\), and all words are read
	from left to right; zero powers and \(\mathrm{id}\) denote
	the empty word.
	
	At stage \(j\), put $N_j^L
	:=
	\varphi_{i_j}(L_{j-1})
	-
	\varepsilon_{i_j}(R_{j-1})$.
	Substitution in \eqref{eq:epsphi1}--\eqref{eq:epsphi2},
	using \eqref{eq:vw-pqmu} and
	\eqref{eq:tau-sigma-eta-identities}, together with the full
	branch words above, gives inductively
	\[
	\begin{array}{c|c|c|c}\hline
		j&N_j^L&\text{left factor}&\text{right factor}\\ \hline
		1&h&
		\mathsf I_1^h&
		\mathsf{II}_1^q\mathsf I_1^u\\ \hline
		2&h&
		\mathsf I_2^h&
		\mathsf{II}_2^{q+d}\mathsf I_2^{u+v}\\ \hline
		3&2h+\eta+q+2d+p&
		\mathsf{II}_1^{2h}
		\mathsf I_1^{\eta+q+2d+p}&
		\mathsf{III}_1^\sigma
		\mathsf{II}_1^{\mu+2u+2v}\\ \hline
		4&h+q+2d+p&
		\mathsf{II}_2^h
		\mathsf I_2^{q+2d+p}&
		\mathsf{II}_2^{\mu+u+v}\\ \hline
		5&h+q+3d+2p+\tau&
		\mathsf{III}_1^h
		\mathsf{II}_1^{q+3d+2p}
		\mathsf I_1^\tau&
		\mathsf{III}_1^{\eta+\mu+u+v}\\ \hline
		6&\mu+p+d&
		\mathsf{II}_2^{p+d}\mathsf I_2^\mu&
		\mathrm{id}\\ \hline
	\end{array}
	\]
	Here \(0\le N_j^L\le\varphi_{i_j}(L_{j-1})\).
	Each lowering move decreases the tensor-factor comparison by one,
	so \eqref{eq:tensor-fe} selects the left factor for the first
	\(N_j^L\) moves and then the right factor.
	Thus each left word is the corresponding prefix of a full
	one-factor string, and each right word is a full string.
	
	At stage~5, let \(L\otimes R\) denote the tensor obtained from
	\(Y^{(4)}\) after the left-factor moves
	\(\mathsf{III}_1^h\mathsf{II}_1^{q+3d+2p}\).
	Then we have
	\[
		(z_2(L),z_3(L))
		=(-q-3d-2p,\sigma),\quad
		\varphi_1(L)=k+\tau+\sigma,\quad
		\varepsilon_1(R)=k+\sigma.
	\]
	If \(\sigma=0\), the next \(\tau\) left-factor moves use
	\(\mathsf I_1\).
	If \(\sigma>0\), then \(\tau=0\), and equality in
	\eqref{eq:tensor-fe} selects the right factor.
	
	Applying the branch words gives
	\begin{equation}\label{eq:classical-lowering-endpoint}
		\begin{aligned}
			Y^{(1)}
			={}&
			(a-h,h,0,0,0,0)
			\otimes
			(0,u+v,\mu,\mu+2q+2d,0,k),
			\\
			Y^{(2)}
			={}&
			(a-h,0,2h,0,0,0)
			\otimes
			(0,0,\mu+2u+2v,\mu,q+d,k),
			\\
			Y^{(3)}
			={}&
			(k+\tau,\mu+q+2d+p-\tau,0,2h,0,0)
			\otimes
			(0,0,0,2\mu+2u+2v,\eta,k+\sigma),
			\\
			Y^{(4)}
			={}&
			(k+\tau,\eta,2q+4d+2p,0,h,0)
			\otimes
			(0,0,0,0,\eta+\mu+u+v,k+\sigma),
			\\
			Y^{(5)}
			={}&
			(k,\mu,q+d,q+2p+3d,0,h)
			\otimes
			\bar{\mathbf1}^{\,r},
			\\
			Y^{(6)}
			={}&
			(k,0,2\mu+q+d,q+d,d+p,h)
			\otimes
			\bar{\mathbf1}^{\,r}
			=
			X^\downarrow.
		\end{aligned}
	\end{equation}
	All listed moves are nonzero, and substitution in
	\eqref{eq:epsphi1}--\eqref{eq:epsphi2} and
	\eqref{eq:tensor-string-functions} gives $\varphi_{i_j}(Y^{(j)})=0$ $(1\le j\le6)$.
	Hence each row gives the maximal lowering string required by
	\eqref{eq:classical-lowering-recursion}.
	
	Finally, we have
	\[
		(\varepsilon_1,\varepsilon_2,\varphi_1,\varphi_2)
		(\bar{\mathbf1}^{\,r})
		=(r,0,0,0),\quad
		\varphi_1(b^\downarrow)
		=k+2\mu+q+d+p=r-u\le r,\quad
		\varphi_2(b^\downarrow)=0.
	\]
	Therefore \eqref{eq:tensor-fe} gives
	\(\widetilde f_1X^\downarrow
	=\widetilde f_2X^\downarrow=0\).
	The displayed path connects \(\hw(X)\) to \(X^\downarrow\),
	so \(X^\downarrow\) belongs to the classical component of \(X\).
	This component is a finite highest-weight \(G_2\)-crystal,
	as in the proof of Proposition~\ref{prop:hw-highest},
	and therefore has a unique lowest element.
\end{proof}

\begin{remark}\label{rem:s-invariance}
	By \eqref{eq:e1}--\eqref{eq:f2}, classical arrows preserve
	\(s\) on each tensor factor.
	Hence, for \(X=b\otimes c\in B_l\otimes B_l\), the parameters
	in \eqref{eq:high-param}--\eqref{eq:r-def} satisfy $a=s(b)$ and $	r=s(c)$.
\end{remark}

\subsection{Defect coordinates and canonical affine string}

\subsubsection{Defect-coordinate identities}

Define the level-independent function
\begin{equation}\label{eq:Theta}
	\Theta:
	\bigcup_{j\ge1}B_j
	\longrightarrow\mathbb Z,
	\quad
	\Theta(b):=M(b)-s(b).
\end{equation}
Superscript \((j)\) denotes color-\(0\) data on \(B_j\).
By \eqref{eq:epsphi0}, we have
\begin{equation}\label{eq:epsphi0-shift}
	\varphi_0^{(j)}(b)=j+\Theta(b)
	\ (j\in\Z_{>0},\ b\in B_j).
\end{equation}

\begin{lemma}\label{lem:affine-drift}
	Let \(j\in\Z_{>0}\) and \(b\in B_j\).
	If \(\widetilde f_0^{(j)}b\ne0\), then $\Theta\bigl(\widetilde f_0^{(j)}b\bigr)
	=
	\Theta(b)-1$.
\end{lemma}

\begin{proof}
	It follows from \textup{(C3)} and \eqref{eq:epsphi0-shift} that $	\Theta\bigl(\widetilde f_0^{(j)}b\bigr)
	=
	\varphi_0^{(j)}
	\bigl(\widetilde f_0^{(j)}b\bigr)-j
	=
	\varphi_0^{(j)}(b)-j-1
	=
	\Theta(b)-1$.
\end{proof}

Fix \(X\in B_l\otimes B_l\), and write
\(U:=U(X)\), \(V:=V(X)\), and \(W:=W(X)\)
as in \eqref{eq:UVW}.
By \eqref{eq:s-bdown}, define
\begin{equation}\label{eq:C-def}
	C_i
	:=
	A_i(b^\downarrow)-s(b^\downarrow)+r
	=
	A_i(b^\downarrow)-a+r
	\ (1\le i\le6).
\end{equation}

\begin{lemma}\label{lem:C-coordinates}
	The coordinates in \eqref{eq:C-def} are
	\begin{equation}\label{eq:C-list}
	\begin{aligned}
		C_1&=V,
		\quad
		C_2=u+\mu-d-k,
		\quad
		C_3=u+\mu-d-k+p-q,\\
		C_4&=u-2\mu-d-k+p-q,
		\quad
		C_5=u+p-k,
		\quad
		C_6=h+u+p-2k.
	\end{aligned}
	\end{equation}
\end{lemma}

\begin{proof}
	At \(b^\downarrow\),
	\eqref{eq:z-def} and \eqref{eq:bdown-def} give $(z_1,z_2,z_3,z_4)
	=
	(h-k,p-q,2\mu+q+d,-\mu)$.
	Moreover, \eqref{eq:r-def} and \eqref{eq:vw-pqmu}--\eqref{eq:h-def}
	 yield $V=r-a=u+\mu-d-h$.
	Substitution in \eqref{eq:A-vector} and
	\eqref{eq:C-def} gives \eqref{eq:C-list}.
\end{proof}

Lemma~\ref{lem:C-coordinates} gives
\begin{equation}\label{eq:C4-C3}
	C_4=C_3-3\mu\le C_3.
\end{equation}

Define
\begin{equation}\label{eq:LQ-def}
	L:=\max\{0,C_2,C_3,C_5\},
	\quad
	Q:=C_6=h+u+p-2k.
\end{equation}

\begin{lemma}
	\label{lem:UVW-compression}
	The following identities hold:
	\begin{align}
		W&=\max\{C_2,C_3,C_5\},
		\quad
		L=\max\{0,W\},
		\label{eq:W-C}\\
		U&=\max\{Q,L-V\}.
		\label{eq:U-QLV}
	\end{align}
\end{lemma}

\begin{proof}
	Since \((p-q)_+=p\) and \(\mu=\min\{v,w\}\),
	\eqref{eq:C-list} gives
	\[
		\max\{C_2,C_3,C_5\}
		=u+p-k+(\mu-d)_+
		=u+v-k-\min\{v,w,d\}
		=W.
	\]
	Equation~\eqref{eq:LQ-def} then gives \(L=\max\{0,W\}\).
	
	By \eqref{eq:r-def}, \eqref{eq:h-def}, and
	\(D^\sharp=2k-2u-v+d\), we obtain
	\[
	a-r-D^\sharp=Q,
	\quad
	Q+k-u+d=u+v-k-V.
	\]
It is easy to see that \(\max\{k-u,k-u-w+d,k-u-v+d\}
=k-u+d-\min\{v,w,d\}\).
	Hence \eqref{eq:M-sharp} and \eqref{eq:UVW} yield
	\[
		U
		=\max\{
		Q,-V,Q+k-u+d-\min\{v,w,d\}
		\}
		=\max\{Q,-V,W-V\}
		=\max\{Q,L-V\}.\qedhere
	\]
\end{proof}

By \eqref{eq:C-def}, \eqref{eq:Theta}, and
\eqref{eq:C-list}--\eqref{eq:W-C}, we have $\Theta(b^\downarrow)+r
=
\max_{1\le i\le6}C_i
=
\max\{V,W,Q\}$.
Hence
\begin{equation}\label{eq:t-defects}
	\max\{0,\Theta(b^\downarrow)+r\}
	=
	\max\{0,V,W,Q\}.
\end{equation}

\subsubsection{The canonical affine string}

For \(X\in B_l\otimes B_l\), use \(b^\downarrow\) and \(r\)
from \eqref{eq:bdown-def} and \eqref{eq:r-def}.
Lemma~\ref{lem:classical-lowest-normal-form} gives
\(b^\downarrow\in B_l\) and \(0\le r\le l\).

For the color-\(0\) segment verifying~\textup{(P1)} in
\cite[Section~7.1]{KMOY}, set
\begin{equation}\label{eq:t-def}
	t:=\bigl(\varphi_0(b^\downarrow)-l+r\bigr)_+.
\end{equation}
Equations~\eqref{eq:epsphi0-shift} and \eqref{eq:t-defects} give
\begin{equation}\label{eq:t-max}
	t
	=
	\bigl(\Theta(b^\downarrow)+r\bigr)_+
	=
	\max\{0,V,W,Q\}.
\end{equation}

For \(n\in\Zge\) with \(n\le t\), set
\(b_n:=\widetilde f_0^{\,n}b^\downarrow\).

\begin{lemma}
	\label{lem:KMOY-zero-schedule}
	The elements \(b_n\) belong to \(B_l\) for \(0\le n\le t\), and
	\begin{align}
		\widetilde f_0^{\,n}
		\bigl(b^\downarrow\otimes\bar{\mathbf1}^{\,r}\bigr)
		&=
		b_n\otimes\bar{\mathbf1}^{\,r},
		\  0\le n\le t,
		\label{eq:zero-schedule-left}\\
		\widetilde f_0^{\,t+j}
		\bigl(b^\downarrow\otimes\bar{\mathbf1}^{\,r}\bigr)
		&=
		b_t\otimes\bar{\mathbf1}^{\,r-j},
		\  0\le j\le r.
		\label{eq:zero-schedule-right}
	\end{align}
\end{lemma}

\begin{proof}
	Since \(0\le r\le l\), \eqref{eq:t-def} gives
	\(0\le t\le\varphi_0(b^\downarrow)\).
	Normality therefore implies \(b_n\in B_l\) for \(0\le n\le t\).
	
	It follows from \eqref{eq:epsphi0} and \eqref{eq:f0} that $	\varepsilon_0(\bar{\mathbf1}^{\,j})=l-j$ $(0\le j\le l)$ and $\widetilde f_0\bar{\mathbf1}^{\,j}=\bar{\mathbf1}^{\,j-1}$ $(1\le j\le l)$.

	For \(0\le n<t\), \textup{(C3)} and \eqref{eq:t-def} give $\varphi_0(b_n)
	=
	\varphi_0(b^\downarrow)-n
	=
	l-r+t-n
	>
	l-r
	=
	\varepsilon_0(\bar{\mathbf1}^{\,r})$.
	Thus \eqref{eq:tensor-fe} selects the left factor for the first
	\(t\) moves, proving \eqref{eq:zero-schedule-left}.
	
	At \(n=t\), we have $\varphi_0(b_t)
	=
	\varphi_0(b^\downarrow)-t
	=
	\min\{\varphi_0(b^\downarrow),l-r\}
	\le l-r$.
	Hence $\varphi_0(b_t)
	\le l-r+j
	=
	\varepsilon_0(\bar{\mathbf1}^{\,r-j})$ $(0\le j<r)$,
	and \eqref{eq:tensor-fe} gives $\widetilde f_0
	\bigl(b_t\otimes\bar{\mathbf1}^{\,r-j}\bigr)
	=
	b_t\otimes\bar{\mathbf1}^{\,r-j-1}$.
	Starting from \eqref{eq:zero-schedule-left} with \(n=t\),
	iterating this identity proves \eqref{eq:zero-schedule-right}.
\end{proof}

Taking \(j=r\) in \eqref{eq:zero-schedule-right} gives
\begin{equation}\label{eq:zero-schedule-endpoint}
	\widetilde f_0^{\,t+r}
	\bigl(b^\downarrow\otimes\bar{\mathbf1}^{\,r}\bigr)
	=
	b_t\otimes\phi.
\end{equation}

For \(0\le n\le t\), define
\begin{equation}\label{eq:Cin-def}
	C_i(n)
	:=
	A_i(b_n)-s(b_n)+r
	\qquad(1\le i\le6).
\end{equation}
Then \(C_i(0)=C_i\).  Set
\begin{equation}\label{eq:mprime-def}
	m':=s(b_t).
\end{equation}

Since \(b_t\in B_l\), we have \(m'\in\{0,\ldots,l\}\).

\begin{lemma}\label{lem:terminal-C1}
	No \(b_n\), \(0\le n\le t\), satisfies \((F_4)\), and $C_1(t)=\min\{0,V+t-L\}$.
\end{lemma}

\begin{proof}
	Write $C(n):=(C_1(n),\ldots,C_6(n))$ and $L_n:=\max\{0,C_2(n),C_3(n),C_5(n)\}$.
	Then \(L_0=L\).
	The entries of \(C(n)\) are integers, and
	\(C_i(n)-C_j(n)=A_i(b_n)-A_j(b_n)\).
	Thus \eqref{eq:F-max-char} identifies \((F_j)\) with the
	first maximal entry \(C_j(n)\).
	For \(0\le n<t\), \eqref{eq:f0} and \eqref{eq:Cin-def} give
	\begin{equation}\label{eq:C-transition-table}
		\begin{array}{c|c}\hline
			\text{condition}&C(n+1)-C(n)\\ \hline
			(F_1)&(-1,-2,-2,-2,-2,-3)\\ \hline
			(F_2)&(0,-1,-2,-2,-1,-2)\\ \hline
			(F_3)&(0,0,-1,-4,-2,-2)\\ \hline
			(F_4)&(0,0,2,-1,-2,-2)\\ \hline
			(F_5)&(0,-1,0,0,-1,-2)\\ \hline
			(F_6)&(1,0,0,0,0,-1)\\ \hline
		\end{array}
	\end{equation}
	
	The table shows that \(C_4-C_3\) is unchanged or decreases
	by \(3\).  Hence \eqref{eq:C4-C3} gives
	\[
	C_4(n)-C_3(n)
	\le C_4(0)-C_3(0)
	=-3\mu\le0
	\ (0\le n\le t).
	\]
	The first-maximum rule therefore excludes \((F_4)\).
	
	If \(t=0\), \eqref{eq:t-max} and \eqref{eq:W-C} give
	\(V\le0\) and \(L=0\), so
	\(C_1(t)=V=\min\{0,V+t-L\}\).
	Assume \(t>0\).
	By \eqref{eq:Cin-def}, \eqref{eq:Theta},
	\eqref{eq:t-max}, and Lemma~\ref{lem:affine-drift}, we have
	\begin{equation}\label{eq:C-max-tn}
		M_n
		:=
		\max_{1\le i\le6}C_i(n)
		=
		\Theta(b_n)+r
		=
		t-n
		\ (0\le n\le t).
	\end{equation}
	In particular, \(M_n\ge1\) for \(n<t\).
	
	Let \(N\) be the first index \(0\le n<t\) for which
	\(b_n\) satisfies \((F_1)\), setting \(N=t\) if no such
	index exists.
	For \(0\le n<N\), only
	\((F_2),(F_3),(F_5),(F_6)\) occur.
	Under \((F_2),(F_3),(F_5)\), the active entry becomes
	\(M_n-1\), while the other entries among \(C_2,C_3,C_5\)
	become at most \(M_n-1\), by integrality and the
	first-maximum rule.
	Since \(M_n\ge1\), this gives
	\(L_n=M_n\) and \(L_{n+1}=M_n-1\).
	An \(F_6\)-move fixes \(L_n\).
	Consequently, \eqref{eq:C-transition-table} yields
	\[
	\bigl(L_{n+1}-L_n,C_1(n+1)-C_1(n)\bigr)
	=
	\begin{cases}
		(-1,0),&(F_2),(F_3),(F_5),\\
		(0,1),&(F_6),
	\end{cases}
	\]
	and therefore $C_1(n)-L_n=V-L+n$ $(0\le n\le N)$.
	
	If \(N=t\), \eqref{eq:C-max-tn} gives \(L_t=0\), hence $C_1(t)=V+t-L\le0$.
	If \(N<t\), then \(C_1(N)=M_N>0\), so \(L_N\le C_1(N)\) and $V+t-L
	=
	C_1(N)-L_N+t-N
	\ge t-N>0$.

	Whenever \((F_1)\) holds at \(b_n\) with \(N\le n<t\),
	\eqref{eq:C-transition-table} gives
	\[
	C_1(n+1)=M_n-1>
	\max_{2\le j\le6}C_j(n+1).
	\]
	Thus \((F_1)\) persists through \(b_t\), and
	\(C_1(t)=M_t=0\).
	In either case,
	\(C_1(t)=\min\{0,V+t-L\}\).
\end{proof}

Since \(A_1(b_t)=0\), equations~\eqref{eq:Cin-def}--\eqref{eq:mprime-def} give \(C_1(t)=r-m'\).
Combining this identity with Lemma~\ref{lem:terminal-C1} yields
\begin{equation}\label{eq:mprime-relation}
	m'-r=\max\{0,L-V-t\}.
\end{equation}

\subsection{Energy evaluation and closed formula}

For \(X\in B_l\otimes B_l\), write \(V:=V(X)\), and use
\(r,t,L,m'\) from \eqref{eq:r-def}, \eqref{eq:t-def},
\eqref{eq:LQ-def}, and \eqref{eq:mprime-def}, respectively.

\begin{proposition}
	\label{prop:canonical-energy}
	For every \(X\in B_l\otimes B_l\), $H_l(X)=t+m'-r=\max\{t,L-V\}$.
\end{proposition}

\begin{proof}
The classical path in Lemma~\ref{lem:classical-lowest-normal-form}
preserves \(H_l\).
Applying \eqref{eq:energy-def} along the affine string in
Lemma~\ref{lem:KMOY-zero-schedule}, with endpoint
\eqref{eq:zero-schedule-endpoint}, gives $H_l(X)
=
H_l\bigl(b^\downarrow\otimes\bar{\mathbf1}^{\,r}\bigr)
=
H_l(b_t\otimes\phi)+t-r$.

	By Remark~\ref{rem:s-invariance} and \eqref{eq:r-def},
	the highest-parameter tuple of \(b_t\otimes\phi\) is
	\((m',0,0,0,0,0)\).
	Equations~\eqref{eq:h-def}--\eqref{eq:bdown-def},
	together with Lemma~\ref{lem:classical-lowest-normal-form},
	therefore identify its classical lowest element as
	\(\bar{\mathbf1}^{\,m'}\otimes\phi\).
	Hence \eqref{eq:energy-def} gives $H_l(b_t\otimes\phi)
	=
	H_l\bigl(\bar{\mathbf1}^{\,m'}\otimes\phi\bigr)$.
	
	For \(1\le j\le m'\), \eqref{eq:epsphi0} gives $\varphi_0(\bar{\mathbf1}^{\,j})
	=
	l+j
	>
	l
	=
	\varepsilon_0(\phi)$.
	Thus \eqref{eq:f0} and \eqref{eq:tensor-fe} give $\widetilde f_0^{\,m'}
	\bigl(\bar{\mathbf1}^{\,m'}\otimes\phi\bigr)
	=
	\phi\otimes\phi$,
	with all \(m'\) moves acting on the left factor.
	By \eqref{eq:energy-def} and \eqref{eq:H-normalization}, we have $H_l\bigl(\bar{\mathbf1}^{\,m'}\otimes\phi\bigr)=m'$.
	Combining these identities with \eqref{eq:mprime-relation},
	we obtain
	\[
	H_l(X)
	=
	t+m'-r
	=
	t+\max\{0,L-V-t\}
	=
	\max\{t,L-V\}.\qedhere
	\]
\end{proof}

\begin{theorem}\label{thm:energy}
	For every \(X\in B_l\otimes B_l\),
	\begin{equation}\label{eq:H-final-three}
		H_l(X)
		=
		\HH(X)
		=
		\max\{U(X),V(X),W(X)\}.
	\end{equation}
\end{theorem}

\begin{proof}
	Proposition~\ref{prop:canonical-energy},
	\eqref{eq:t-max}, and
	Lemma~\ref{lem:UVW-compression} give
	\[
	\begin{aligned}
		H_l(X)
		=\max\{t,L-V\}
		=\max\{0,V,W,Q,L-V\}
		=\max\{0,U,V,W\}.
	\end{aligned}
	\]
	Since \(L\ge0\) and \(U\ge L-V\), we have \(U\ge -V\).
	Thus \(\max\{U,V\}\ge |V|\ge0\), so the zero branch
	is redundant.
	Equation~\eqref{eq:H-three} now yields
	\eqref{eq:H-final-three}.
\end{proof}

Proposition~\ref{prop:nine-branch} gives the nine-branch formula
\eqref{eq:H-nine} for \(H_l\).
In particular, \(H_l(X)\ge |a-r|\ge0\).

\begin{corollary}\label{cor:level-stability}
	Let \(j\in\Z_{>0}\) with \(j\le l\). Then $H_l\big|_{B_j\otimes B_j}=H_j$
	under the coordinate inclusion \(B_j\subseteq B_l\).
\end{corollary}

\begin{proof}
	The classical operator formulas and string functions are
	level-independent.
	By Remark~\ref{rem:s-invariance}, the reduction
	\eqref{eq:HW-recursion} starting in \(B_j\otimes B_j\)
	remains there and agrees at levels \(j\) and \(l\).
	It therefore gives the same highest parameters and the same
	\(U,V,W\) at both levels.
	Equation~\eqref{eq:H-final-three} proves the claim.
\end{proof}

\subsection{The direct twelve-coordinate formula}

Evaluating \eqref{eq:HW-recursion} in the original coordinates
eliminates the highest parameters in \eqref{eq:H-nine}; by
Theorem~\ref{thm:energy}, this gives \(H_l\) in the twelve
coordinates of \(b\otimes c\).

Let $b=x:=(x_1,x_2,x_3,\bar x_3,\bar x_2,\bar x_1)$ and $c=y:=(y_1,y_2,y_3,\bar y_3,\bar y_2,\bar y_1)$
be elements of \(B_l\).
With \(s\) as in \eqref{eq:s-def}, set
\begin{equation}\label{eq:twelve-sx}
	s_x:=s(b),
	\quad
	s_y:=s(c).
\end{equation}

Define
\begin{align*}
	\rho_x
	&:=
	\bigl(
	x_2-x_3-(\bar x_2-\bar x_3)_+
	\bigr)_+,
	\quad
	\sigma_x
	:=
	(\bar x_3-\bar x_2)_+,
	\quad
	\tau_x
	:=
	\frac{
		\bigl(x_3-\bar x_3+2\sigma_x\bigr)_+
	}{2},
	\\
	p_x
	&:=
	x_2-\rho_x+\tau_x,
	\quad
	q_x
	:=
	x_3+\sigma_x-2\tau_x.
\end{align*}
The left-factor thresholds for the first five raising stages,
with colors \(1,2,1,2,1\), are
\begin{align*}
	\kappa_1
	&:=
	x_1+
	\bigl(
	x_3-x_2+(\bar x_2-\bar x_3)_+
	\bigr)_+,
	\quad
	\kappa_2
	:=
	x_2-\rho_x+
	\frac{
		\bigl(\bar x_3-x_3-2\sigma_x\bigr)_+
	}{2},
	\\
	\kappa_3
	&:=
	x_1+\rho_x+(q_x-p_x)_+,
	\quad
	\kappa_4
	:=
	\min\{p_x,q_x\},
	\quad
	\kappa_5
	:=
	x_1+\rho_x+(p_x-q_x)_+.
\end{align*}
The sixth-stage threshold is zero.
All these quantities are nonnegative, since
\[
0\le\rho_x\le x_2,
\quad
q_x
=
\min\{x_3+\sigma_x,\bar x_3-\sigma_x\}
\ge0.
\]

For the first raising stage, define
\begin{align*}
	\sigma_1
	&:=
	(\bar y_3-\bar y_2)_+,
	\quad
	\rho_1
	:=
	\bigl(
	y_2-y_3-(\bar y_2-\bar y_3)_+
	\bigr)_+,
	\quad
	r_{13}
	:=
	(\bar y_1-\kappa_1)_+,
	\\
	r_{12}
	&:=
	(\bar y_1+\sigma_1-\kappa_1)_+
	-
	(\bar y_1-\kappa_1)_+,
	\quad
	r_{11}
	:=
	(\bar y_1+\sigma_1+\rho_1-\kappa_1)_+
	-
	(\bar y_1+\sigma_1-\kappa_1)_+.
\end{align*}

For the second stage, define
\begin{align*}
	q_{22}
	&:=
	(\bar y_2+r_{13}-\kappa_2)_+,
	\\
	q_{21}
	&:=
	\left(
	\bar y_2+r_{13}
	+
	\frac{
		\bigl(y_3-\bar y_3+2r_{12}\bigr)_+
	}{2}
	-\kappa_2
	\right)_+
	-
	(\bar y_2+r_{13}-\kappa_2)_+.
\end{align*}

For the third stage, define
\begin{align*}
	\sigma_3
	&:=
	\bigl(
	\bar y_3-\bar y_2-r_{12}-r_{13}+3q_{22}
	\bigr)_+,
	\\
	\rho_3
	&:=
	\Bigl(
	y_2-y_3-r_{11}-r_{12}+3q_{21}
	-
	\bigl(
	\bar y_2-\bar y_3+r_{12}+r_{13}-3q_{22}
	\bigr)_+
	\Bigr)_+,
\end{align*}
\begin{align*}
	r_{33}
	&:=
	(\bar y_1-r_{13}-\kappa_3)_+,
	\quad
	r_{32}
	:=
	(\bar y_1-r_{13}+\sigma_3-\kappa_3)_+
	-
	(\bar y_1-r_{13}-\kappa_3)_+,
	\\
	r_{31}
	&:=
	(\bar y_1-r_{13}+\sigma_3+\rho_3-\kappa_3)_+
	-
	(\bar y_1-r_{13}+\sigma_3-\kappa_3)_+.
\end{align*}

For the fourth stage, put
\begin{align*}
	\omega_4
	&:=
	\bar y_2+r_{13}-q_{22}+r_{33},
	\\
	\Delta_4
	&:=
	y_3-\bar y_3
	+2r_{12}
	-2q_{21}-2q_{22}
	+2r_{32},
\end{align*}
and define
\begin{align*}
	q_{42}
	&:=
	(\omega_4-\kappa_4)_+,
	\\
	q_{41}
	&:=
	\left(
	\omega_4+\frac{(\Delta_4)_+}{2}-\kappa_4
	\right)_+
	-
	(\omega_4-\kappa_4)_+.
\end{align*}

For the fifth stage, define
\begin{align*}
	\sigma_5
	&:=
	\Bigl(
	\bar y_3-\bar y_2
	-r_{12}-r_{13}-r_{32}-r_{33}
	+3q_{22}+3q_{42}
	\Bigr)_+,
	\\
	\rho_5
	&:=
	\Bigl(
	y_2-y_3
	-r_{11}-r_{12}-r_{31}-r_{32}
	+3q_{21}+3q_{41}
	\\
	&\hspace{1cm}
	-
	\bigl(
	\bar y_2-\bar y_3
	+r_{12}+r_{13}+r_{32}+r_{33}
	-3q_{22}-3q_{42}
	\bigr)_+
	\Bigr)_+,
\end{align*}
\begin{align*}
	r_{53}
	&:=
	(\bar y_1-r_{13}-r_{33}-\kappa_5)_+,
	\\
	r_{52}
	&:=
	(\bar y_1-r_{13}-r_{33}+\sigma_5-\kappa_5)_+
	-
	(\bar y_1-r_{13}-r_{33}-\kappa_5)_+,
	\\
	r_{51}
	&:=
	(\bar y_1-r_{13}-r_{33}
	+\sigma_5+\rho_5-\kappa_5)_+
	-
	(\bar y_1-r_{13}-r_{33}
	+\sigma_5-\kappa_5)_+.
\end{align*}

For the sixth stage, define
\begin{align*}
	q_{61}
	&:=
	\frac12
	\Bigl(
	y_3-\bar y_3
	+2r_{12}
	-2q_{21}-2q_{22}
	+2r_{32}
	-2q_{41}-2q_{42}
	+2r_{52}
	\Bigr)_+,
	\\
	q_{62}
	&:=
	\bar y_2+r_{13}-q_{22}
	+r_{33}-q_{42}+r_{53}.
\end{align*}

The positive-part formulas and monotonicity of
\(t\mapsto(t)_+\) give \(r_{ij}\ge0\) and
\(q_{21},q_{22},q_{41},q_{42},q_{61}\ge0\).
Moreover,
\[
\omega_4
=
\min\{\bar y_2+r_{13},\kappa_2\}+r_{33}
\ge0,
\quad
q_{62}
=
\min\{\omega_4,\kappa_4\}+r_{53}
\ge0.
\]
Set
\begin{align*}
	R_1
	&:=
	r_{11}+r_{31}+r_{51},
	\quad
	R_2
	:=
	r_{12}+r_{32}+r_{52},
	\\
	R_3
	&:=
	r_{13}+r_{33}+r_{53},
	\quad
	Q_{\mathrm{coord}}
	:=
	q_{21}+q_{41}+q_{61}.
\end{align*}
The definition of \(q_{62}\) gives
\begin{equation}\label{eq:twelve-q2-total}
	q_{22}+q_{42}+q_{62}
	=
	\bar y_2+R_3.
\end{equation}

Define
\begin{align}
	\mathcal F_1
	&:=
	s_x+y_1-\bar y_3-2\bar y_2-3\bar y_1
	+R_1+R_2+R_3,
	\label{eq:twelve-F1}\\
	\mathcal F_2
	&:=
	s_x-\bar y_3-2\bar y_2-2\bar y_1+R_2,
	\notag\\
	\mathcal F_3
	&:=
	s_x-\frac{3y_3+\bar y_3}{2}
	-\bar y_2-2\bar y_1
	+3Q_{\mathrm{coord}}-R_2+R_3,
	\notag\\
	\mathcal F_4
	&:=
	s_x-y_2-\frac{y_3+\bar y_3}{2}
	-\bar y_2-2\bar y_1
	+R_1+R_3,
	\notag\\
	\mathcal F_5
	&:=
	s_x-s_y,
	\notag\\
	\mathcal F_6
	&:=
	s_y-s_x,
	\notag\\
	\mathcal F_7
	&:=
	y_1-\bar y_1+R_1+R_3,
	\notag\\
	\mathcal F_8
	&:=
	y_1+y_2-y_3-\bar y_1
	+3Q_{\mathrm{coord}}-R_2+R_3,
	\notag\\
	\mathcal F_9
	&:=
	y_1+y_2+\frac{y_3-\bar y_3}{2}
	-\bar y_2-\bar y_1+R_2.
	\label{eq:twelve-F9}
\end{align}

The parity conditions in \eqref{eq:Bl-set} give
\begin{equation}\label{eq:twelve-parity}
	x_3\pm\bar x_3\in2\mathbb Z,
	\quad
	y_3\pm\bar y_3\in2\mathbb Z,
	\quad
	3y_3+\bar y_3\in2\mathbb Z.
\end{equation}
Positive parts preserve \(\mathbb Z\) and \(2\mathbb Z\).
Every numerator divided by \(2\) is an even expression in
\eqref{eq:twelve-parity}, or its negative, plus twice an
already integral quantity, possibly enclosed in a positive part.
Induction through the definitions of \(r_{j\nu}\), \(q_{j\nu}\), and \(\mathcal F_j\) therefore gives
\[
r_{j\nu}\in\Zge
\ (j\in\{1,3,5\},\ 1\le\nu\le3),
\quad
q_{j\nu}\in\Zge
\ (j\in\{2,4,6\},\ 1\le\nu\le2),
\]
and $\mathcal F_j\in\mathbb Z$ $(1\le j\le9)$.

By the definitions from \eqref{eq:twelve-sx} through
\eqref{eq:twelve-F9}, the functions
\(\mathcal F_1,\ldots,\mathcal F_9\) are piecewise-linear on
\(\mathbb R^6\times\mathbb R^6\).
Define
\begin{equation}\label{eq:direct-coordinate-H}
	\mathscr H_{\mathrm{coord}}:
	\mathbb R^6\times\mathbb R^6
	\longrightarrow\mathbb R,
	\qquad
	\mathscr H_{\mathrm{coord}}(x,y)
	:=
	\max_{1\le j\le9}\mathcal F_j.
\end{equation}

\begin{theorem}
	\label{thm:twelve-coordinate-energy}
	For \(b,c\in B_l\), with coordinate vectors \(x,y\),
	\begin{equation}\label{eq:twelve-coordinate-energy}
		H_l(b\otimes c)
		=
		\mathscr H_{\mathrm{coord}}(x,y).
	\end{equation}
\end{theorem}

\begin{proof}
	Apply \eqref{eq:HW-recursion} to \(X=b\otimes c\), and write
	\(X^{(j)}=b^{(j)}\otimes c^{(j)}\) for \(0\le j\le6\).
	At the \(j\)-th raising stage, Lemma~\ref{lem:Pi-emax} gives
	\(\varepsilon_{i_j}(b^{(j-1)})\) and $\bigl(
	\varepsilon_{i_j}(c^{(j-1)})
	-
	\varphi_{i_j}(b^{(j-1)})
	\bigr)_+$
	raising moves on the first and second factors, respectively.
	
	Put
	\(\lambda_x:=\bar x_3-\sigma_x+2\bar x_2+2\bar x_1\).
	Then we have
	\[
		\lambda_x
		=\min\{\bar x_3,\bar x_2\}+2\bar x_2+2\bar x_1
		\in\Zge,\quad
		q_x+\lambda_x
		=x_3+\bar x_3+2\bar x_2+2\bar x_1-2\tau_x
		\in2\Zge.
	\]
	Equations~\eqref{eq:E1m} and \eqref{eq:E2m} give
	\[
	\begin{array}{c|c|c}\hline
		j&b^{(j)}
		&\varphi_{i_{j+1}}(b^{(j)})\\ \hline
		0&
		(x_1,x_2,x_3,\bar x_3,\bar x_2,\bar x_1)
		&\kappa_1\\ \hline
		1&
		(x_1+\rho_x,x_2-\rho_x,x_3+\sigma_x,
		\bar x_3-\sigma_x,\bar x_2+\bar x_1,0)
		&\kappa_2\\ \hline
		2&
		(x_1+\rho_x,p_x,q_x,\lambda_x,0,0)
		&\kappa_3\\ \hline
		3&
		(\kappa_5,\kappa_4,q_x+\lambda_x,0,0,0)
		&\kappa_4\\ \hline
		4&
		(\kappa_5,s_x-\kappa_5,0,0,0,0)
		&\kappa_5\\ \hline
		5&
		(s_x,0,0,0,0,0)
		&0\\ \hline
	\end{array}
	\]
	The sixth operation fixes the first factor, so
	\(b^{(6)}=b^{(5)}\).
	
	For a one-factor raising string with ordered block capacities
	\(A,B,C\in\Zge\) and left-factor threshold \(\kappa\in\Zge\),
	the right-factor counts are
	\[
	(A+B+C-\kappa)_+-(B+C-\kappa)_+,\quad
	(B+C-\kappa)_+-(C-\kappa)_+,\quad
	(C-\kappa)_+.
	\]
	For two ordered blocks of capacities \(A,B\), the corresponding
	counts are
	\[
	(A+B-\kappa)_+-(B-\kappa)_+,
	\quad
	(B-\kappa)_+.
	\]
	Induction on \(j\), using \eqref{eq:E1m}--\eqref{eq:E2m}, gives
	\[
	\begin{array}{c|c|c|c}\hline
		j&i_j&
		\text{ordered capacities of }c^{(j-1)}
		&\text{right-factor moves}\\ \hline
		1&1&
		(\rho_1,\sigma_1,\bar y_1)
		&(r_{11},r_{12},r_{13})\\ \hline
		2&2&
		\left(
		\dfrac{(y_3-\bar y_3+2r_{12})_+}{2},
		\bar y_2+r_{13}
		\right)
		&(q_{21},q_{22})\\ \hline
		3&1&
		(\rho_3,\sigma_3,\bar y_1-r_{13})
		&(r_{31},r_{32},r_{33})\\ \hline
		4&2&
		\left(
		\dfrac{(\Delta_4)_+}{2},
		\omega_4
		\right)
		&(q_{41},q_{42})\\ \hline
		5&1&
		(\rho_5,\sigma_5,\bar y_1-r_{13}-r_{33})
		&(r_{51},r_{52},r_{53})\\ \hline
		6&2&
		(q_{61},q_{62})
		&(q_{61},q_{62})\\ \hline
	\end{array}
	\]
	The thresholds in these six rows are
	\(\kappa_1,\ldots,\kappa_5,0\), respectively.
	
	Consequently, we have
	\[
	c^{(j)}-c^{(j-1)}
	=
	\begin{cases}
		(r_{j1},-r_{j1},r_{j2},-r_{j2},r_{j3},-r_{j3}),
		&j\in\{1,3,5\},\\[1mm]
		(0,q_{j1},-2q_{j1},2q_{j2},-q_{j2},0),
		&j\in\{2,4,6\}.
	\end{cases}
	\]
	Summing these six changes and using
	\eqref{eq:twelve-q2-total} gives
	\[
	c^{(6)}
	=
	\bigl(
	y_1+R_1,\,
	y_2-R_1+Q_{\mathrm{coord}},\,
	y_3+R_2-2Q_{\mathrm{coord}},\,
	\bar y_3-R_2+2\bar y_2+2R_3,\,
	0,\,
	\bar y_1-R_3
	\bigr).
	\]
	
	Comparing \(X^{(6)}=\hw(b\otimes c)\) with
	\eqref{eq:high-param} and using Remark~\ref{rem:s-invariance}
	gives
	\[
	\begin{aligned}
		a&=s_x,
		\quad
		u=y_1+R_1,
		\quad
		v=y_2-R_1+Q_{\mathrm{coord}},
		\quad
		w=y_3+R_2-2Q_{\mathrm{coord}},
		\\
		k&=\bar y_1-R_3,
		\quad
		r=s_y,
		\quad
		d=
		\frac{\bar y_3-y_3}{2}
		+\bar y_2+R_3-R_2+Q_{\mathrm{coord}}.
	\end{aligned}
	\]
	Substituting these parameters into \eqref{eq:H-nine}
	identifies its nine branches with
	\(\mathcal F_1,\ldots,\mathcal F_9\), in the same order.
	Theorem~\ref{thm:energy} and \eqref{eq:direct-coordinate-H}
	therefore give \eqref{eq:twelve-coordinate-energy}.
\end{proof}

In \eqref{eq:twelve-coordinate-energy}, \(l\) enters only through
\(b,c\in B_l\).
For \(l=1\), the set \(B_1\) consists of the following eight elements:
\[
\begin{aligned}
	v_a&=(1,0,0,0,0,0),\quad
	v_b=(0,1,0,0,0,0),\quad
	v_c=(0,0,2,0,0,0),\quad
	v_d=(0,0,1,1,0,0),\\
	v_e&=(0,0,0,2,0,0),\quad
	v_f=(0,0,0,0,1,0),\quad
	v_g=(0,0,0,0,0,1),\quad
	v_h=\phi=(0,0,0,0,0,0).
\end{aligned}
\]

\begin{remark}
	Equation~\eqref{eq:twelve-coordinate-energy} gives the following
	table, with rows indexing the left tensor factor and columns
	the right:
	\begin{equation}\label{eq:H1}
		\begin{array}{c|c|c|c|c|c|c|c|c}\hline
			H_1&v_a&v_b&v_c&v_d&v_e&v_f&v_g&v_h\\ \hline
			v_a&2&1&1&0&0&0&0&1\\ \hline
			v_b&2&2&1&1&1&0&0&1\\ \hline
			v_c&2&2&2&1&1&1&0&1\\ \hline
			v_d&2&2&2&1&1&1&0&1\\ \hline
			v_e&2&2&2&2&2&1&1&1\\ \hline
			v_f&2&2&2&2&2&2&1&1\\ \hline
			v_g&2&2&2&2&2&2&2&1\\ \hline
			v_h&1&1&1&1&1&1&1&0\\ \hline
		\end{array}
	\end{equation}
	The matrix agrees entrywise with \cite[Figure~6]{Dombos}.
	After simultaneously reordering the rows and columns as
	\((v_h,v_a,\ldots,v_g)\), it agrees with \cite[(3.3)]{FHKS}.
\end{remark}

\section{\texorpdfstring{Similarity transfer to type $G_2^{(1)}$}
	{Similarity transfer to type G2(1)}}
\label{sec:G21-similarity-transfer}

Throughout this section, fix \(l\in\Z_{>0}\).
For \(n\in\Z_{>0}\), write \(B_n^D:=B_n\) and \(H_n^D:=H_n\).
We use superscripts \(D,G\) for crystal operators and string
functions, retaining \(\mathcal A,M,D\) and
\(\mathcal A_G,M_G,D_G\) for the respective coordinate functions.
The similarity factors \(m_i\) are given by
\eqref{eq:G-similarity-factors}.

\subsection{Similarity and local-energy transfer}

\subsubsection{One-factor similarity}

Define
\begin{equation}\label{eq:G-similarity-map}
	v_l:B_l^G\longrightarrow B_{3l}^D,
	\quad
	v_l(b):=3b,
\end{equation}
where multiplication is coordinatewise, and extend \(v_l\) by
\(v_l(0):=0\).

For \(b=(x_1,x_2,x_3,\bar x_3,\bar x_2,\bar x_1)\in B_l^G\),
\eqref{eq:G-Bl-set} gives
\[
3b\in\Zge^6,
\quad
3x_3\equiv3\bar x_3\!\!\!\!\pmod2,
\quad
s(3b)=3s(b)\le3l.
\]
Hence \(v_l\) is a well-defined injection into \(B_{3l}^D\).
Linearity of \(s,\mathcal A,D\) and positive homogeneity of
the maximum give
\begin{equation}\label{eq:G-coordinate-homogeneity}
	\begin{aligned}
		s(v_l(b))=3s(b),
		\quad
		\mathcal A(v_l(b))=3\mathcal A_G(b),
		\quad
		M(v_l(b))=3M_G(b),
		\quad
		D(v_l(b))=3D_G(b).
	\end{aligned}
\end{equation}
All type-\(D_4^{(3)}\) operators and string functions at
\(v_l(b)\) are taken at level \(3l\).

\begin{proposition}
	\label{prop:G-one-factor-similarity}
	For \(b\in B_l^G\) and \(i\in I\),
	\begin{equation}\label{eq:operator-identities}
	v_l(\widetilde e_i^G b)
	=
	(\widetilde e_i^D)^{m_i}v_l(b),
	\quad
	v_l(\widetilde f_i^G b)
	=
	(\widetilde f_i^D)^{m_i}v_l(b),
	\end{equation}
	\begin{equation}\label{eq:G-one-factor-string-scaling}
		\varepsilon_i^D(v_l(b))
		=
		m_i\varepsilon_i^G(b),
		\quad
		\varphi_i^D(v_l(b))
		=
		m_i\varphi_i^G(b).
	\end{equation}
\end{proposition}

\begin{proof}
	Let \(\widehat B_{3l}\) be the \(1/3\)-scaled
	type-\(D_4^{(3)}\) crystal of
	\cite[Section~2, equations~(2.1)--(2.5)]{MMO},
	with operators \(\widehat e_i,\widehat f_i\).
	Dilation \(b\mapsto3b\), extended by \(0\mapsto0\),
	identifies \(\widehat B_{3l}\) with \(B_{3l}^D\)
	and conjugates these operators to
	\(\widetilde e_i^D,\widetilde f_i^D\), respectively.
	
	By \eqref{eq:G-Bl-set} and \cite[(3.1)]{MMO},
	\(B_l^G=V_l\subset\widehat B_{3l}\).
	Moreover, \cite[Proposition~1(2)]{MMO} gives $\widetilde e_i^G b=(\widehat e_i)^{m_i}b$ and $\widetilde f_i^G b=(\widehat f_i)^{m_i}b$.
	Conjugating these identities by dilation yields
	\eqref{eq:operator-identities}.
	
	Taking the type-\(D_4^{(3)}\) level to be \(3l\),
	we obtain from \eqref{eq:epsphi0}, \eqref{eq:epsphi1},
	\eqref{eq:epsphi2}, \eqref{eq:G-epsphi2}, and
	\eqref{eq:G-coordinate-homogeneity} that
	\[
	\bigl(
	\varepsilon_i^D(v_l(b)),
	\varphi_i^D(v_l(b))
	\bigr)
	=
	\begin{cases}
		3\bigl(\varepsilon_i^G(b),\varphi_i^G(b)\bigr),
		&i=0,1,\\
		\bigl(\varepsilon_2^G(b),\varphi_2^G(b)\bigr),
		&i=2.
	\end{cases}
	\]
	Equation~\eqref{eq:G-similarity-factors} now gives
	\eqref{eq:G-one-factor-string-scaling}.
\end{proof}

\subsubsection{Tensor-product compatibility}

Define
\[
v_l^{\otimes2}:
B_l^G\otimes B_l^G
\longrightarrow
B_{3l}^D\otimes B_{3l}^D,
\quad
v_l^{\otimes2}(b\otimes c)
:=
v_l(b)\otimes v_l(c)
\]
and set \(v_l^{\otimes2}(0):=0\).

\begin{lemma}
	\label{lem:G-no-switching}
	Let \(X=b\otimes c\in B_l^G\otimes B_l^G\) and \(i\in I\).
	
	If \(\widetilde f_i^G X\ne0\), then
	\begin{equation}\label{eq:G-tensor-f-similarity}
	v_l^{\otimes2}(\widetilde f_i^G X)
	=
	(\widetilde f_i^D)^{m_i}v_l^{\otimes2}(X),
	\end{equation}
	and all \(m_i\) type-\(D_4^{(3)}\) moves act on the tensor factor
	selected by \(\widetilde f_i^G\).
	
	If \(\widetilde e_i^G X\ne0\), then
	\begin{equation}\label{eq:G-tensor-e-similarity}
		v_l^{\otimes2}(\widetilde e_i^G X)
		=
		(\widetilde e_i^D)^{m_i}v_l^{\otimes2}(X),
	\end{equation}
	and all \(m_i\) type-\(D_4^{(3)}\) moves act on the tensor factor
	selected by \(\widetilde e_i^G\).
\end{lemma}

\begin{proof}
	Set
	\(\Delta:=\varphi_i^G(b)-\varepsilon_i^G(c)\in\mathbb Z\).
	By \eqref{eq:G-one-factor-string-scaling}, the initial
	type-\(D\) comparison is
	\[
	\varphi_i^D(v_l(b))
	-
	\varepsilon_i^D(v_l(c))
	=
	m_i\Delta.
	\]
	
	Under either nonzero hypothesis, the operator on the selected
	factor is nonzero.
	Proposition~\ref{prop:G-one-factor-similarity} therefore
	guarantees \(m_i\) nonzero type-\(D\) one-factor moves.
	
	By \textup{(C3)}, each nonzero lowering move decreases the
	comparison by \(1\), so after \(k\) such moves it is
	\(m_i\Delta-k\).
	The type-\(G\) lowering rule selects the left factor exactly
	when \(\Delta\ge1\).
	For \(0\le k<m_i\),
	\[
	\begin{array}{ll}
		m_i\Delta-k\ge m_i-k>0,
		&\Delta\ge1,\\
		m_i\Delta-k\le-k\le0,
		&\Delta\le0.
	\end{array}
	\]
	Thus \eqref{eq:tensor-fe} shows inductively that all \(m_i\)
	lowering moves are nonzero and select the same factor as
	\(\widetilde f_i^G\).
	
	By \textup{(C2)}, each nonzero raising move increases the
	comparison by \(1\), so after \(k\) such moves it is
	\(m_i\Delta+k\).
	The type-\(G\) raising rule selects the left factor exactly
	when \(\Delta\ge0\).
	For \(0\le k<m_i\),
	\[
	\begin{array}{ll}
		m_i\Delta+k\ge k\ge0,
		&\Delta\ge0,\\
		m_i\Delta+k\le-m_i+k<0,
		&\Delta\le-1.
	\end{array}
	\]
	The same induction shows that all \(m_i\) raising moves are
	nonzero and select the same factor as \(\widetilde e_i^G\).
	In particular, when \(\Delta=0\), lowering acts on the right
	and raising on the left, as required by \eqref{eq:tensor-fe}.
	
In each case, all \(m_i\) type-\(D\) moves act on the image
of the factor selected by the corresponding type-\(G\) operator, so Proposition~\ref{prop:G-one-factor-similarity} yields
\eqref{eq:G-tensor-f-similarity} and
\eqref{eq:G-tensor-e-similarity}.
\end{proof}

\begin{corollary}
	\label{cor:G-tensor-string-scaling}
For \(X\in B_l^G\otimes B_l^G\) and \(i\in I\),
\(\varepsilon_i^D\bigl(v_l^{\otimes2}(X)\bigr)
=m_i\varepsilon_i^G(X)\) and
\(\varphi_i^D\bigl(v_l^{\otimes2}(X)\bigr)
=m_i\varphi_i^G(X)\).
\end{corollary}

\begin{proof}
	Write \(X=b\otimes c\).
	By \eqref{eq:tensor-string-functions} and \textup{(C1)},
	for either type, we have
	\[
		\varepsilon_i(b\otimes c)
		=
		\varepsilon_i(b)
		+
		\bigl(\varepsilon_i(c)-\varphi_i(b)\bigr)_+,\quad
		\varphi_i(b\otimes c)
		=
		\varphi_i(c)
		+
		\bigl(\varphi_i(b)-\varepsilon_i(c)\bigr)_+.
	\]
	Equation~\eqref{eq:G-one-factor-string-scaling} and
	positive homogeneity of \(t\mapsto(t)_+\) prove both identities.
\end{proof}

By Corollary~\ref{cor:G-tensor-string-scaling} and normality,
the operator identities in Lemma~\ref{lem:G-no-switching}
extend to the formal null cases.
Thus, for every \(X\in B_l^G\otimes B_l^G\) and \(i\in I\),
\[
v_l^{\otimes2}(\widetilde e_i^G X)
=
(\widetilde e_i^D)^{m_i}v_l^{\otimes2}(X),
\quad
v_l^{\otimes2}(\widetilde f_i^G X)
=
(\widetilde f_i^D)^{m_i}v_l^{\otimes2}(X).
\]

\subsubsection{Transfer of the normalized local energy}

The zero six-tuple \(\phi\) belongs to \(B_l^G\) and satisfies
\(v_l(\phi)=\phi\).
Define
\begin{equation}\label{eq:G-energy-transfer}
	H_l^G:
	B_l^G\otimes B_l^G
	\longrightarrow
	\tfrac13\mathbb Z,
	\quad
	H_l^G(X)
	:=
	\frac13
	H_{3l}^D\bigl(v_l^{\otimes2}(X)\bigr).
\end{equation}

\begin{theorem}
	\label{thm:G-energy-transfer}
	The function \(H_l^G\) is integer-valued and is the unique
	local energy function on \(B_l^G\otimes B_l^G\) satisfying
	\(H_l^G(\phi\otimes\phi)=0\).
\end{theorem}

\begin{proof}
	Let \(X=b\otimes c\), \(i\in I\), and
	\(Y:=\widetilde f_i^G X\ne0\).
	By Lemma~\ref{lem:G-no-switching}, \(v_l^{\otimes2}\)
	sends this edge to \(m_i\) nonzero type-\(D_4^{(3)}\) edges,
	all on the factor selected by \(\widetilde f_i^G\).
	Hence \eqref{eq:G-energy-transfer}, \eqref{eq:energy-def},
	\eqref{eq:tensor-fe}, and \(m_0=3\) give
	\begin{equation}\label{eq:G-energy-edge-increments}
		H_l^G(Y)-H_l^G(X)
		=
		\begin{cases}
			0,
			&i=1,2,\\
			-1,
			&i=0,\quad
			\varphi_0^G(b)>\varepsilon_0^G(c),\\
			1,
			&i=0,\quad
			\varphi_0^G(b)\le\varepsilon_0^G(c).
		\end{cases}
	\end{equation}
	
	By \eqref{eq:H-normalization} at level \(3l\) and
	\(v_l(\phi)=\phi\), we have \(H_l^G(\phi\otimes\phi)=0\).
	Perfectness of \(B_l^G\) and \textup{(P1)} imply that
	\(B_l^G\otimes B_l^G\) is connected.
	Summing the integer increments in
	\eqref{eq:G-energy-edge-increments}, with appropriate signs,
	along an unoriented path from \(\phi\otimes\phi\) to \(X\)
	gives \(H_l^G(X)\in\mathbb Z\) for every \(X\).
	Thus \(H_l^G\) is a normalized local energy function.
	Uniqueness follows from
	Proposition~\ref{prop:perfect-energy-existence}.
\end{proof}

\subsection{\texorpdfstring{Closed local-energy formulas in type $G_2^{(1)}$}
		{Closed local-energy formulas in type G2(1)}}

\subsubsection{Highest-parameter formula}

For \(i\in\{1,2\}\) and \(Y\in B_l^G\otimes B_l^G\), set
\((\widetilde e_i^G)^{\max}Y
:=(\widetilde e_i^G)^{\varepsilon_i^G(Y)}Y\).

For \(X\in B_l^G\otimes B_l^G\), set \(X_G^{(0)}:=X\) and define
\begin{equation}\label{eq:G-HW-recursion}
	X_G^{(k)}
	:=
	(\widetilde e_{i_k}^G)^{\max}X_G^{(k-1)}
	\ (1\le k\le6),
\end{equation}
where \((i_1,\ldots,i_6)=(1,2,1,2,1,2)\).  Set
\begin{equation}\label{eq:G-hw-def}
	\hw_G(X):=X_G^{(6)}.
\end{equation}
Write \(\hw_D:=\hw\) for the type-\(D_4^{(3)}\) highest reduction on
\(B_{3l}^D\otimes B_{3l}^D\) defined by
\eqref{eq:HW-recursion} and \eqref{eq:hw-def}.

\begin{lemma}
	\label{lem:G-highest-synchronization}
	For every \(X\in B_l^G\otimes B_l^G\), $\hw_D\bigl(v_l^{\otimes2}(X)\bigr)
	=
	v_l^{\otimes2}\bigl(\hw_G(X)\bigr)$.
\end{lemma}

\begin{proof}
	Let \(Y\in B_l^G\otimes B_l^G\), \(i\in\{1,2\}\), and
	put \(N:=\varepsilon_i^G(Y)\).
	By normality,
	\((\widetilde e_i^G)^{n+1}Y\ne0\) for \(0\le n<N\).
	Iterating \eqref{eq:G-tensor-e-similarity} and applying
	Corollary~\ref{cor:G-tensor-string-scaling} give
	\[
	\begin{aligned}
		v_l^{\otimes2}
		\bigl((\widetilde e_i^G)^{\max}Y\bigr)
		=
		(\widetilde e_i^D)^{m_iN}v_l^{\otimes2}(Y)
		=
		(\widetilde e_i^D)^{
			\varepsilon_i^D(v_l^{\otimes2}(Y))}
		v_l^{\otimes2}(Y).
	\end{aligned}
	\]
	Iterating this identity along \eqref{eq:G-HW-recursion},
	and using Lemma~\ref{lem:Pi-emax},
	\eqref{eq:HW-recursion}--\eqref{eq:hw-def}, and
	\eqref{eq:G-hw-def}, proves the claim.
\end{proof}

By \eqref{eq:G-HW-recursion}--\eqref{eq:G-hw-def},
\(\hw_G(X)\) belongs to the classical component of \(X\).
Lemma~\ref{lem:G-highest-synchronization},
Proposition~\ref{prop:hw-highest}, and
Corollary~\ref{cor:G-tensor-string-scaling} give
\(\varepsilon_i^G(\hw_G(X))=0\) for \(i=1,2\).

For any \(G_2\)-highest element \(Z\) in this component,
Lemma~\ref{lem:G-no-switching} places
\(v_l^{\otimes2}(Z)\) in the classical component of
\(v_l^{\otimes2}(X)\), and
Corollary~\ref{cor:G-tensor-string-scaling} gives
\(\varepsilon_i^D(v_l^{\otimes2}(Z))=0\) for \(i=1,2\).
Proposition~\ref{prop:hw-highest} and
Lemma~\ref{lem:G-highest-synchronization} therefore yield $v_l^{\otimes2}(Z)
=
\hw_D\bigl(v_l^{\otimes2}(X)\bigr)
=
v_l^{\otimes2}\bigl(\hw_G(X)\bigr)$.
The injectivity of \(v_l^{\otimes2}\) gives \(Z=\hw_G(X)\),
proving uniqueness.

\begin{lemma}
	\label{lem:G-highest-shape}
	An element \(Y\in B_l^G\otimes B_l^G\) is annihilated by
	\(\widetilde e_1^G\) and \(\widetilde e_2^G\) if and only if there
	are unique \(a,u,v,w,d,k\) such that
	\begin{equation}\label{eq:G-highest-shape}
		Y
		=
		(a,0,0,0,0,0)
		\otimes
		(u,v,w,w+2d,0,k),
	\end{equation}
	where
	\begin{equation}\label{eq:G-highest-lattice}
		a,u,k\in\Zge,
		\quad
		v,w,d\in\tfrac13\Zge,
		\quad
		v-w,\ w+2d\in\mathbb Z,
	\end{equation}
	and $a\le l$, $u+v+w+d+k\le l$, $k+w+2d+(v-w)_+\le a$.
\end{lemma}

\begin{proof}
	Write \(Y=b\otimes c\).
	Normality and \eqref{eq:tensor-string-functions} give
	\[
	\widetilde e_i^G Y=0
	\iff
	\varepsilon_i^G(b)=0
	\ \text{and}\
	\varepsilon_i^G(c)\le\varphi_i^G(b)
	\ (i=1,2).
	\]
	By \eqref{eq:epsphi1} and \eqref{eq:G-epsphi2},
	\[
	\varepsilon_1^G(b)=\varepsilon_2^G(b)=0
	\iff
	b=(a,0,0,0,0,0)
	\ (a\in\Zge),
	\]
	in which case
	\(\varphi_1^G(b)=a\) and \(\varphi_2^G(b)=0\).
	Thus the remaining string conditions are
	\(\varepsilon_2^G(c)=0\) and \(\varepsilon_1^G(c)\le a\).
	
	By \eqref{eq:G-epsphi2}, \(\varepsilon_2^G(c)=0\) gives
	\(c=(u,v,w,w+2d,0,k)\) with \(u,v,w,d,k\ge0\).
	The coordinate conditions in \eqref{eq:G-Bl-set} give $u,k\in\Zge$, $v,w\in\tfrac13\Zge$ and $v-w,\ w+2d\in\mathbb Z$,
	while the parity condition in \eqref{eq:G-Bl-set} gives
	\(6d=3(w+2d)-3w\in2\mathbb Z\), hence
	\(d\in\tfrac13\Zge\).
	Therefore \eqref{eq:G-highest-lattice} holds.
	Moreover, \eqref{eq:epsphi1} and \eqref{eq:s-def} give
	\[
	\varepsilon_1^G(c)=k+w+2d+(v-w)_+,
	\quad
	s(b)=a,
	\quad
	s(c)=u+v+w+d+k.
	\]
	These yield the three stated inequalities.
	
	Conversely, \eqref{eq:G-highest-lattice} gives nonnegative
	coordinates in \(\tfrac13\mathbb Z\) and the integrality
	conditions in \eqref{eq:G-Bl-set}; the parity condition follows from
	\(3(w+2d)-3w=6d\in2\mathbb Z\).
	The level bounds place both factors in \(B_l^G\).
	The string formulas and the last inequality give $\varepsilon_i^G(b)=0$ and $\varepsilon_i^G(c)\le\varphi_i^G(b)$ for $i=1, 2$.
	The tensor criterion therefore gives
	\(\widetilde e_1^G Y=\widetilde e_2^G Y=0\).
	Uniqueness follows directly from the coordinates in
	\eqref{eq:G-highest-shape}.
\end{proof}

For \(X\in B_l^G\otimes B_l^G\), let \(a,u,v,w,d,k\) be the
parameters given by Lemma~\ref{lem:G-highest-shape} for
\(Y=\hw_G(X)\), and set \(r:=u+v+w+d+k\).
The lattice conditions \eqref{eq:G-highest-lattice} and the
level bound in Lemma~\ref{lem:G-highest-shape} give
\(r=u+k+(v-w)+2(w+2d)-3d\in\Zge\) and \(r\le l\).

By \eqref{eq:z-def}--\eqref{eq:A-vector}, and
\(k-u-w-2d\le k-u\), set
\begin{equation}\label{eq:G-UVW}
	\begin{aligned}
		M_G^\sharp
		&:=
		M_G(u,v,w,w+2d,0,k)\\
		&=
		\max\{
		0,\,
		k-u,\,
		k-u-w+d,\,
		k-u-v+d,\,
		2k-2u-v+d
		\},
		\\
		U_G
		&:=
		a+M_G^\sharp-r-(2k-2u-v+d),
		\\
		V_G
		&:=
		r-a,
		\\
		W_G
		&:=
		r-\bigl(2k+w+d+\min\{v,w,d\}\bigr).
	\end{aligned}
\end{equation}

\begin{theorem}
	\label{thm:G-closed-energy}
	For every \(X\in B_l^G\otimes B_l^G\), $H_l^G(X)=\max\{U_G,V_G,W_G\}$.
\end{theorem}

\begin{proof}
	By Lemma~\ref{lem:G-highest-synchronization}, we have
	\[
	\hw_D\bigl(v_l^{\otimes2}(X)\bigr)
	=
	(3a,0,0,0,0,0)
	\otimes
	(3u,3v,3w,3w+6d,0,3k).
	\]
	All six highest parameters are therefore multiplied by \(3\).
	By positive homogeneity, \eqref{eq:UVW} and \eqref{eq:G-UVW} give
	\[
	U\bigl(v_l^{\otimes2}(X)\bigr)=3U_G,\quad
	V\bigl(v_l^{\otimes2}(X)\bigr)=3V_G,\quad
	W\bigl(v_l^{\otimes2}(X)\bigr)=3W_G.
	\]
	Equations~\eqref{eq:G-energy-transfer} and
	\eqref{eq:H-final-three} therefore give
	\[
		H_l^G(X)
		=
		\frac13
		H_{3l}^D\bigl(v_l^{\otimes2}(X)\bigr)
		=
		\frac13\max\{3U_G,3V_G,3W_G\}
		=
		\max\{U_G,V_G,W_G\}.\qedhere
	\]
\end{proof}

By Theorem~\ref{thm:G-closed-energy}, expanding
\eqref{eq:G-UVW} gives
\begin{equation}\label{eq:G-nine-branch-energy}
	\begin{aligned}
		H_l^G(X)=\max\bigl\{&
		a+u-w-2d-3k,\,
		a-w-2d-2k,\,
		a-2w-d-2k,\\
		&
		a-v-w-d-2k,\,
		a-r,\,
		r-a,\,
		u-k,\,
		u+v-w-k,\,
		u+v-d-k
		\bigr\}.
	\end{aligned}
\end{equation}

Every branch in \eqref{eq:G-nine-branch-energy} is integral.
Indeed, \eqref{eq:G-highest-lattice} gives
\(a,u,k,v-w,w+2d,3d\in\mathbb Z\), while
\[
\begin{aligned}
	2w+d
	&=2(w+2d)-3d\in\mathbb Z,\\
	v+w+d
	&=(v-w)+(2w+d)\in\mathbb Z,\\
	v-d
	&=(v-w)+(w+2d)-3d\in\mathbb Z.
\end{aligned}
\]
Together with \(r\in\mathbb Z\), these identities prove the claim.

\begin{corollary}
	For every \(X\in B_l^G\otimes B_l^G\),
	\(H_l^G(X)\ge0\).
\end{corollary}

\begin{proof}
	Equation~\eqref{eq:G-nine-branch-energy} gives
	\(H_l^G(X)\ge\max\{a-r,r-a\}=|a-r|\ge0\).
\end{proof}

\subsubsection{Direct twelve-coordinate formula}

For \(b,c\in B_l^G\), write
\[
\begin{aligned}
	b=x:=(x_1,x_2,x_3,\bar x_3,\bar x_2,\bar x_1),\quad
	c=y:=(y_1,y_2,y_3,\bar y_3,\bar y_2,\bar y_1).
\end{aligned}
\]
The piecewise-linear functions
\eqref{eq:direct-coordinate-H} and
\eqref{eq:twelve-sx}--\eqref{eq:twelve-F9} on
\(\mathbb R^6\times\mathbb R^6\) retain their notation on
\(B_l^G\times B_l^G\).

\begin{theorem}
	\label{thm:G-direct-coordinate-energy}
	For every \(b,c\in B_l^G\),
	\begin{equation}\label{eq:G-direct-coordinate-energy}
		H_l^G(b\otimes c)
		=
		\mathscr H_{\mathrm{coord}}(x,y)
		=
		\max_{1\le j\le9}\mathcal F_j(x,y).
	\end{equation}
\end{theorem}

\begin{proof}
	The formulas \eqref{eq:twelve-sx}--\eqref{eq:direct-coordinate-H}
	involve only linear combinations, positive parts, minima, and maxima.
	They therefore give $\mathscr H_{\mathrm{coord}}(\lambda x,\lambda y)
	=
	\lambda\mathscr H_{\mathrm{coord}}(x,y)$ $(\lambda\ge0)$.
Since \(3b,3c\in B_{3l}^D\) by \eqref{eq:G-similarity-map}, we have
	\[
	\begin{aligned}
		H_l^G(b\otimes c)
		&=
		\frac13H_{3l}^D(3b\otimes3c)
		&&\text{by \eqref{eq:G-energy-transfer}}\\
		&=
		\frac13\mathscr H_{\mathrm{coord}}(3x,3y)
		&&\text{by Theorem~\ref{thm:twelve-coordinate-energy}}\\
		&=
		\mathscr H_{\mathrm{coord}}(x,y).
	\end{aligned}
	\]
	The equality with
	\(\max_{1\le j\le9}\mathcal F_j(x,y)\)
	is \eqref{eq:direct-coordinate-H}.
\end{proof}

The formula depends on \(l\) only through \(b,c\in B_l^G\).
Its auxiliary quantities may be fractional, but
\(\mathscr H_{\mathrm{coord}}(x,y)\) is integral by
Theorem~\ref{thm:G-energy-transfer}.

\begin{corollary}
	\label{cor:G-level-stability}
	Let \(j\in\Z_{>0}\) with \(j\le l\).
	Then \(H_l^G\big|_{B_j^G\otimes B_j^G}=H_j^G\)
	under the coordinate inclusion \(B_j^G\subseteq B_l^G\).
\end{corollary}

\begin{proof}
	For \(b,c\in B_j^G\),
	Theorem~\ref{thm:G-direct-coordinate-energy}
	applied at levels \(j\) and \(l\) gives
	\(H_l^G(b\otimes c)
	=\mathscr H_{\mathrm{coord}}(x,y)
	=H_j^G(b\otimes c)\).
\end{proof}

\begin{remark}
	Theorems~\ref{thm:twelve-coordinate-energy} and
	\ref{thm:G-direct-coordinate-energy} express the two normalized
	energies by the same function
	\(\mathscr H_{\mathrm{coord}}\) on the coordinate domains
	\eqref{eq:Bl-set} and \eqref{eq:G-Bl-set}.
	The coordinate inclusion \(B_l^D\subseteq B_l^G\) is strict,
	since
	\((0,0,0,\tfrac23,\tfrac23,0)\in B_l^G\setminus B_l^D\).
	Thus the full energy matrices at the same level have different sizes.
	
	By \eqref{eq:G-similarity-map} and
	\eqref{eq:G-energy-transfer}, after reindexing by \(v_l\),
	the type-\(G_2^{(1)}\) energy matrix is one third of the
	principal submatrix of the level-\(3l\) type-\(D_4^{(3)}\)
	energy matrix indexed by \(v_l(B_l^G)\):
	\begin{equation}\label{eq:G-energy-matrix-transfer}
		\bigl(
		H_l^G(b\otimes c)
		\bigr)_{b,c\in B_l^G}
		=
		\frac13
		\bigl(
		H_{3l}^D\bigl(v_l(b)\otimes v_l(c)\bigr)
		\bigr)_{b,c\in B_l^G}.
	\end{equation}
\end{remark}

\begin{example}
	Enumeration of \eqref{eq:Bl-set} and
	\eqref{eq:G-Bl-set} gives
	\(\lvert B_3^D\rvert=112\) and \(\lvert B_1^G\rvert=15\).
	Order \(B_3^D=\{b_1^D,\ldots,b_{112}^D\}\) lexicographically by
	\((x_1,x_2,x_3,\bar x_3,\bar x_2,\bar x_1)\), and
	\(B_1^G=\{b_1^G,\ldots,b_{15}^G\}\) so that
	\(v_1(b_1^G),\ldots,v_1(b_{15}^G)\) increase lexicographically.
	Define \(\nu=(\nu_1,\ldots,\nu_{15})\) by
	\(v_1(b_i^G)=b_{\nu_i}^D\) for \(1\le i\le15\).
	The similarity images and their indices are
	\[
	\begin{array}{c|c|c}
		i&v_1(b_i^G)&\nu_i\\ \hline
		1 &(0,0,0,0,0,0)&1\\
		2 &(0,0,0,0,0,3)&4\\
		3 &(0,0,0,0,3,0)&10\\
		4 &(0,0,0,2,2,0)&16\\
		5 &(0,0,0,4,1,0)&19\\
		6 &(0,0,0,6,0,0)&20\\
		7 &(0,0,3,1,1,0)&43\\
		8 &(0,0,3,3,0,0)&44\\
		9 &(0,0,6,0,0,0)&50\\
		10&(0,1,1,1,1,0)&63\\
		11&(0,1,1,3,0,0)&64\\
		12&(0,1,4,0,0,0)&70\\
		13&(0,2,2,0,0,0)&76\\
		14&(0,3,0,0,0,0)&77\\
		15&(3,0,0,0,0,0)&112
	\end{array}.
	\]

	Set
	\(\mathsf H_1^G
	:=
	\bigl(H_1^G(b_i^G\otimes b_j^G)\bigr)_{1\le i,j\le15}\).
	Direct evaluation of \eqref{eq:G-direct-coordinate-energy} gives
	\begingroup
	\scriptsize
	\setlength{\arraycolsep}{2.5pt}
	\[
	\mathsf H_1^G
	=
	\begin{pmatrix}
		0&1&1&1&1&1&1&1&1&1&1&1&1&1&1\\
		1&2&2&2&2&2&2&2&2&2&2&2&2&2&2\\
		1&1&2&2&2&2&2&2&2&2&2&2&2&2&2\\
		1&1&1&1&1&2&1&2&2&1&2&2&2&2&2\\
		1&1&1&1&1&2&1&2&2&1&2&2&2&2&2\\
		1&1&1&1&1&2&1&2&2&1&2&2&2&2&2\\
		1&0&1&1&1&1&1&1&2&1&1&2&2&2&2\\
		1&0&1&1&1&1&1&1&2&1&1&2&2&2&2\\
		1&0&1&1&1&1&1&1&2&1&1&2&2&2&2\\
		1&0&0&0&0&1&0&1&1&0&1&1&1&2&2\\
		1&0&0&0&0&1&0&1&1&0&1&1&1&2&2\\
		1&0&0&0&0&1&0&1&1&0&1&1&1&2&2\\
		1&0&0&0&0&1&0&1&1&0&1&1&1&2&2\\
		1&0&0&0&0&1&0&1&1&0&1&1&1&2&2\\
		1&0&0&0&0&0&0&0&1&0&0&1&1&1&2
	\end{pmatrix}.
	\]
	\endgroup

After simultaneously reordering the rows and columns in the order 
\[
(b_1^G,b_{15}^G,b_{14}^G,b_{13}^G,b_{12}^G,b_{11}^G,
b_9^G,b_8^G,b_{10}^G,b_6^G,b_7^G,b_5^G,b_4^G,b_3^G,b_2^G),
\]
the resulting matrix agrees entrywise with the energy matrix given in
\cite[(3.5)]{FHKS}.
	
	Let $\mathsf H_3^D
	:=
	\bigl(
	H_3^D(b_p^D\otimes b_q^D)
	\bigr)_{1\le p,q\le112}$ and $\mathsf H_3^D[\nu,\nu]
	:=
	\bigl(
	H_3^D(b_{\nu_i}^D\otimes b_{\nu_j}^D)
	\bigr)_{1\le i,j\le15}$.
\(\mathsf H_3^D[\nu,\nu]\) is the principal submatrix ordered by
\(v_1(b_1^G),\ldots,v_1(b_{15}^G)\).
By \eqref{eq:G-energy-matrix-transfer} at \(l=1\),
	\(\mathsf H_1^G=\frac13\mathsf H_3^D[\nu,\nu]\).
	
	For a sample entry, take
	\(b_4^G=(0,0,0,\frac23,\frac23,0)\) and
	\(b_6^G=(0,0,0,2,0,0)\).  Since \(v_1(b_4^G)=b_{16}^D\) and
	\(v_1(b_6^G)=b_{20}^D\), equations~\eqref{eq:G-HW-recursion}--\eqref{eq:G-hw-def},
	together with Lemma~\ref{lem:G-highest-synchronization}, give
	\[
	\begin{aligned}
		\hw_G(b_4^G\otimes b_6^G)
		&=
		(1,0,0,0,0,0)\otimes(1,0,0,0,0,0),\\
		\hw_D(b_{16}^D\otimes b_{20}^D)
		&=
		(3,0,0,0,0,0)\otimes(3,0,0,0,0,0).
	\end{aligned}
	\]
	The corresponding parameter tuples
	\((a,u,v,w,d,k,r)\) are
	\((1,1,0,0,0,0,1)\) and \((3,3,0,0,0,0,3)\), respectively.
	Hence \eqref{eq:G-nine-branch-energy} and \eqref{eq:H-nine} give
	\[
		H_1^G(b_4^G\otimes b_6^G)
		=\max\{2,1,1,1,0,0,1,1,1\}=2,\ \
		H_3^D(b_{16}^D\otimes b_{20}^D)
		=\max\{6,3,3,3,0,0,3,3,3\}=6.
	\]
	Consequently, $H_1^G(b_4^G\otimes b_6^G)
	=
	\frac13
	H_3^D\bigl(
	v_1(b_4^G)\otimes v_1(b_6^G)
	\bigr)$.
\end{example}


\begin{thebibliography}{99}
	\bibitem{BFKL}
	G. Benkart, I. Frenkel, S.-J. Kang, H. Lee,
	\emph{Level 1 perfect crystals and path realizations of basic representations at $q=0$},
	Int. Math. Res. Not. \textbf{2006} (2006), 10312.
	
	\bibitem{Dombos}
	B. Dombos,
	\emph{A new grounded partition identity of type $D_4^{(3)}$},
	SIGMA \textbf{22} (2026), 042.
	
	\bibitem{FHKS}
	Z. Fan, S. Han, S.-J. Kang, Y.-S. Shin,
	\emph{Young wall construction of level-1 highest weight crystals over
		$U_q(D_4^{(3)})$ and $U_q(G_2^{(1)})$},
	J. Algebra \textbf{655} (2024), 376--404.
	
	\bibitem{HanClassical}
	S. Han,
	\emph{Closed formulas for energy functions on tensor squares of higher-level perfect crystals in classical affine types},
	arXiv:2607.20934  (2026).
	
	\bibitem{JayneMisra}
	R. L. Jayne, K. C. Misra,
	\emph{On Demazure crystals for $U_q(G_2^{(1)})$},
	Proc. Amer. Math. Soc. \textbf{139} (2011), 2343--2356.
	
	\bibitem{KMN1}
	S.-J. Kang, M. Kashiwara, K. C. Misra, T. Miwa, T. Nakashima,
	A. Nakayashiki,
	\emph{Affine crystals and vertex models},
	Int. J. Mod. Phys. A \textbf{7} (1992), 449--484.
	
	\bibitem{KMN2}
	S.-J. Kang, M. Kashiwara, K. C. Misra, T. Miwa, T. Nakashima,
	A. Nakayashiki,
	\emph{Perfect crystals of quantum affine Lie algebras},
	Duke Math. J. \textbf{68} (1992), 499--607.
	
	\bibitem{KashiwaraSimilarity}
	M. Kashiwara,
	\emph{Similarity of crystal bases},
	Lie Algebras and Their Representations, Contemp. Math. \textbf{194},
	American Mathematical Society, Providence, RI, 1996, 177--186.
	
	\bibitem{KMOY}
	M. Kashiwara, K. C. Misra, M. Okado, D. Yamada,
	\emph{Perfect crystals for $U_q(D_4^{(3)})$},
	J. Algebra \textbf{317} (2007), 392--423.
	
	\bibitem{LittelmannRO}
	P. Littelmann,
	\emph{Paths and root operators in representation theory},
	Ann. of Math. (2) \textbf{142} (1995), 499--525.
	
	\bibitem{MMO}
	K. C. Misra, M. Mohamad, M. Okado,
	\emph{Zero action on perfect crystals for $U_q(G_2^{(1)})$},
	SIGMA \textbf{6} (2010), 022.
	
	\bibitem{MOW}
	K. C. Misra, M. Okado, E. A. Wilson,
	\emph{Soliton cellular automaton associated with $G_2^{(1)}$ crystal base},
	J. Math. Phys. \textbf{53} (2012), 013510.
	
	\bibitem{Scrimshaw2016}
	T. Scrimshaw,
	\emph{A crystal to rigged configuration bijection and the filling map
		for type $D_4^{(3)}$},
	J. Algebra \textbf{448} (2016), 294--349.
	
	\bibitem{YamadaSCA}
	D. Yamada,
	\emph{Scattering rule in soliton cellular automaton associated with crystal base of $U_q(D_4^{(3)})$},
	J. Math. Phys. \textbf{48} (2007), 043509.
	
	\bibitem{Yamane}
	S. Yamane,
	\emph{Perfect crystals of $U_q(G_2^{(1)})$},
	J. Algebra \textbf{210} (1998), 440--486.
\end{thebibliography}
\end{document}